\documentclass[11pt]{article}

\usepackage{graphicx}
\usepackage{latexsym,amsmath,amsfonts,amscd, amsthm, dsfont}
\usepackage{bm,color}
\usepackage{epsfig,verbatim,epstopdf,graphics}
\usepackage{subfigure}
\usepackage{changebar}
\usepackage{multirow}

\usepackage[ruled,vlined]{algorithm2e}
\SetKwComment{Comment}{$\triangleright$\ }{} 
\DontPrintSemicolon

\usepackage{yhmath}
 \usepackage{booktabs} 
 \usepackage{tikz}
\usepackage{verbatim}
\usetikzlibrary{arrows,backgrounds,snakes,shapes}
 \numberwithin{equation}{section}

\graphicspath{{./}{./figure/}}
\allowdisplaybreaks

\newtheoremstyle{plainNoItalics}{}{}{\normalfont}{}{\bfseries}{.}{ }{}

\theoremstyle{plain}
\newtheorem{thm}{Theorem}[section]

\theoremstyle{plainNoItalics}

\newtheorem{rem}[thm]{Remark}

\newtheorem{exa}[thm]{Example}

\newcommand{\bV}{{\bf V}}

\newcommand{\Dt}{{\Delta t }}

\newcommand{\be}{\begin{eqnarray}}
\newcommand{\ee}{\end{eqnarray}}
\newcommand{\beno}{\begin{eqnarray*}}
\newcommand{\eeno}{\end{eqnarray*}}

\makeatletter

\newcommand{\Rmnum}[1]{\expandafter\@slowromancap\romannumeral #1@}
\makeatother

\begin{document}

\baselineskip=1.8pc



\begin{center}
{\bf
{High Order Semi-Lagrangian Discontinuous Galerkin Method Coupled with Runge-Kutta Exponential Integrators for Nonlinear Vlasov Dynamics}
}
\end{center}

\vspace{.2in}
\centerline{
Xiaofeng Cai\footnote{
 Department of Mathematical Sciences, University of Delaware, Newark, DE, 19716. E-mail: xfcai@udel.edu.
},
Sebastiano Boscarino\footnote{
Department of Mathematics and Computer Science, University of Catania, Catania, 95127, E-mail: boscarino@dmi.unict.it.
},
Jing-Mei Qiu\footnote{Corresponding Author. Department of Mathematical Sciences, University of Delaware, Newark, DE, 19716. E-mail: jingqiu@udel.edu. Research of first and last author is supported by NSF grant NSF-DMS-1818924, Air Force Office of Scientific Computing FA9550-18-1-0257 and University of Delaware.}
}

\bigskip
\noindent
{\bf Abstract.}
In this paper, we propose a semi-Lagrangian discontinuous Galerkin method coupled with Runge-Kutta exponential integrators (SLDG-RKEI) for nonlinear Vlasov dynamics. The commutator-free Runge-Kutta (RK) exponential integrators (EI) were proposed by Celledoni, et al. (FGCS, 2003). In the nonlinear transport setting, the RKEI can be used to decompose the evolution of the nonlinear transport into a composition of a sequence of linearized dynamics. The resulting linearized transport equations can be solved by the semi-Lagrangian (SL) discontinuous Galerkin (DG) method proposed in Cai, et al. (JSC, 2017). The proposed method can achieve high order spatial accuracy via the SLDG framework, and high order temporal accuracy via the RK EI. Due to the SL nature, the proposed SLDG-RKEI method is not subject to the CFL condition, thus they have the potential in using larger time-stepping sizes than those in the Eulerian approach. Inheriting advantages from the SLDG method, the proposed SLDG-RKEI schemes are mass conservative, positivity-preserving, have no dimensional splitting error, perform well in resolving complex solution structures, and can be evolved with adaptive time stepping sizes.
We show the performance of the SLDG-RKEI algorithm by classical test problems for the nonlinear Vlasov-Poisson system, as well as the Guiding center Vlasov model. Though that it is not our focus of this paper to explore the SLDG-RKEI scheme for nonlinear hyperbolic conservation laws that develop shocks, we show some preliminary results on schemes' performance on the Burgers' equation.

\vfill

{\bf Key Words:} Semi-Lagrangian; Discontinuous Galerkin; Runge-Kutta exponential integrators; Vlasov-Poisson;  Guiding center Vlasov model; Mass conservative; Positivity-preserving; Adaptive time-stepping algorithm.
\newpage

\section{Introduction}

In this paper, we consider the following two nonlinear Vlasov models. The first is the nonlinear Vlasov-Poisson system
\begin{equation}
f_t + vf_x + E(x,t) f_v =0,
\label{vlasov}
\end{equation}
\begin{equation}
E(x,t) = -\phi_x, \
-\phi_{xx}(x,t) = \rho(x,t).
\label{poisson}
\end{equation}
Here $x$ and $v$ are the coordinates in the phase space $(x,v)\in\Omega_x \times \mathbb{R}$; the electron distribution function $f(x,v,t)$ is the probability distribution function describing the probability of finding a particle with velocity $v$ at position $x$ and at time $t$.
The electric field $E=-\phi_x$, where the self-consistent electrostatic potential $\phi$ is determined by Poisson's equation \eqref{poisson}. $\rho(x,t) = \int_{ \mathbb{R} } f(x,v,t) dv -1$ denotes charge density, with the assumption that infinitely massive ions are uniformly distributed in the background.
The second model is the guiding center Vlasov model which describes a highly magnetized plasma in the transverse plane of a tokamak \cite{shoucri1981two,crouseilles2009conservative}:
\begin{equation}
\rho_t + \nabla \cdot(  \mathbf{E}^{\bot} \rho ) = 0,
\label{guiding}
\end{equation}
\begin{equation}
-\Delta \Phi = \rho,\
\mathbf{E}^{\bot} =  ( -\Phi_y,\Phi_x ),
\label{poisson2d}
\end{equation}
where the unknown variable $\rho$ denotes the charge density of the plasma, and the electric field $\mathbf{E}$ depends on $\rho$ via the Poisson equation. Both models can be written in the form of
\begin{equation}
u_t + \nabla \cdot( \mathbf{P}( u;\mathbf{x},t ) u ) = 0, \ (\mathbf{x},t)\in\mathbb{R}^d\times[0,T],
\label{general_nonlinear}
\end{equation}
where  $u:\mathbb{R}^d\times[0,T]\rightarrow  \mathbb{R}$, $\mathbf{P}( u;\mathbf{x},t  ) = ( P_1(u;\mathbf{x},t ),\cdots,P_d(u;\mathbf{x},t ) )^T$ with $P_i: \mathbb{R}\times\mathbb{R}^d\times[0,T]\rightarrow \mathbb{R},i=1,\cdots,d$ are velocity fields that are $u$-dependent for nonlinear dynamics.

Popular mesh-based approaches for the above mentioned nonlinear transport models are the Eulerian and semi-Lagrangian (SL) approaches. Eulerian methods are usually built by a spatial discretization coupled by a temporal time discretization of the partial differential equations via the method-of-lines approach; while SL methods are designed taking into account characteristics tracing. Despite the complication in building in the characteristics tracking mechanism, when properly designed, SL methods can circumvent the stringent $CFL$ constraint in an Eulerian approach, thus achieve computational savings by taking larger time-stepping sizes. In particular,  when performing a time integration of a system, one would take $\Delta t$ to be $\min (\Delta t_{acc}, \Delta t_{stab})$, where $\Delta t_{acc}$ stands for the time stepping size  from accuracy consideration, while $\Delta t_{stab}$ is the time stepping size from stability consideration. By working with the SL approach, $\Delta t_{stab}$ could be greatly relaxed, so that one can take the time stepping size with accuracy consideration only, leading to computational savings by taking a larger time stepping size. 
The computational savings, result from a larger time stepping size, could become even more significant for a nonlinear model such as the guiding-center Vlasov model, where the dominant computational cost in a time step is the elliptic solver for the Poisson equation.

The semi-Lagrangian approach is becoming attractive in many application domains including the plasmas simulations \cite{sonnendrucker2004vlasov,besse2008wavelet,crouseilles2014new,grandgirard20165d,kormann2019massively}, climate modeling \cite{lin1997explicit,lauritzen2010conservative}, fluid mechanics \cite{xiu2001semi,restelli2006semi,celledoni2016high,bonaventura2018fully,peixoto2019semi}, and kinetic modelling \cite{dimarco2015multiscale,groppi2016boundary}.
Besides the capability of large time step size, the scheme should also be high order and  mass conservative for satisfying the demand for applications.
Recently, a high order conservative SL finite difference WENO scheme is proposed in \cite{qiu_shu_sl}; it can be applied for high dimensional problems via the dimensional splitting method.
However, in the frame of the finite difference scheme, it is highly nontrivial to propose a high order conservative SL schemes without dimensional splitting error.
Thus finite volume schemes or discontinuous Galerkin (DG) schemes are used as the spatial discretization of the semi-Lagrangian method since mass conservation can be assured conveniently. It is well known that the Runge-Kutta (RK) discontinuous Galerkin method (DG) \cite{cockburn2001runge} is popular for the problems \eqref{general_nonlinear} because of its low numerical dissipation, compactness, flexibility for boundary and parallel implementation, superconvergence, high resolution for discontinuities. Then semi-Lagrangian DG schemes are proposed for achieving the goal of overcoming the CFL condition constraint and inheriting as many good properties from DG as possible. There are two classes of SLDG methods: one is base on the weak Galerkin form from the upstream cells for which the mass conservation comes from proper trace of upstream elements \cite{rossmanith2011positivity, cai2016high}, the other is base on the flux function, where the flux functions are evaluated following characteristics \cite{restelli2006semi, qiu2011positivity} and the mass conservation follows naturally from the unique definition of fluxes at element interfaces. The proposed method in this paper belongs to the first class. Recently, many SLDG schemes are proposed based on a  dimensional splitting strategy (see \cite{qiu2011positivity,crouseilles2011discontinuous,rossmanith2011positivity,kometa2011semi,Guo2013discontinuous,besse2017adaptive,einkemmer2019performance} and the references therein).
For problems (e.g. the guiding center Vlasov model), the dimensional splitting error of these schemes may dominate the computational errors  \cite{cai2019comparison}.
Because of this, the authors use a non-splitting SLDG algorithm in \cite{cai2016high}  for its mass conservation, up to third order spatial accuracy,  compactness, non-oscillatory as well as positivity preserving (PP); we refer the reader to \cite{restelli2006semi,tumolo2012semi,lee2016high,bosler2019conservative} for the references of other nonsplitting SLDG schemes.
The formulation of the SLDG scheme in \cite{cai2016high} is similar in spirit to the characteristic Galerkin weak formulation in the ELLAM where the treatment of general boundary conditions are given  \cite{celia1990eulerian,wang1999ellam,wang2006eulerian,cheng2010preliminary}. For convergence and error analysis, the optimal convergence and superconvergence of SLDG schemes for linear convection equations in one space dimension are shown in \cite{yang2020optimal}; the convergence of high order numerical schemes is discussed in \cite{einkemmer2014convergence} for nonlinear Vlasov dynamics. Theoretical study of the proposed method for nonlinear dynamics will be the subject of our future investigation.

For the use of high order SL methods for nonlinear dynamic such as \eqref{general_nonlinear},
the tracking of characteristics with high order temporal accuracy is still a nontrivial issue.
Research efforts have been made to accurately track characteristics for nonlinear Vlasov dynamics \cite{qiu2017high}; yet the problem-dependent procedures become more and more complicated when higher order temporal accuracy is desired. On the other hand, a class of commutator-free Runge-Kutta (RK) exponential integrator (EI) are proposed in the context of SL schemes to solve the nonlinear convection-dominated problems \cite{celledoni2003commutator}. The RKEI framework constructs schemes
by decomposing the nonlinear dynamic process into a sequence of linearized linear solvers; and the high order temporal accuracy is achieved by matching order conditions. The RKEI schemes are represented in the form of Butcher tableaus; as such, the schemes can be implemented in a black-box manner, as  in implementing the RK time discretization in an Eulerian approach. 
In this paper, we propose to apply the SLDG algorithm for linear transport problems \cite{cai2016high} to couple with the RKEI for nonlinear Vlasov dynamics.
When the velocity field in a nonlinear problem is being linearized around the DG solution at that time step, the velocity field becomes discontinuous with jumps at element edges. If issues from discontinuous velocity field arise, one possible remedy is to apply smoothness-increasing accuracy-conserving filters as in \cite{li2016smoothness, cockburn2003enhanced} to DG solutions for the velocity field, which is one of our future directions to further pursue.  Finally, we would  like to mention an alternative approach in treating the nonlinearity in characteristics tracing. In \cite{huang2012eulerian,huang2016semi}, characteristics are approximated to a first order, where the error of such approximation is being taken into account by a correction term by the flux function.

The rest of this paper is organized as follows. In Section \ref{section:method}, we propose to couple the SLDG and RKEI method; in Section \ref{sec3}, the performance of the proposed method is shown through extensive numerical tests. Finally, concluding remarks are made in Section \ref{section:conclusion}.

\section{SLDG-RKEI schemes}
\label{section:method}

In this section, we will present the proposed SLDG-RKEI method that combines the SLDG schemes \cite{cai2016high} with the high order RKEI in \cite{celledoni2003commutator} for solving nonlinear transport problems.
We will first review the SLDG scheme for linear transport problems.
To extend this SLDG solver for nonlinear transport problems,
we start by illustrating a first order SLDG-RKEI scheme that updates the solution by using the SLDG to solve a linearized transport equation.
In order to achieve high order accuracy in time,
the SLDG is coupled with the high order RKEI that decomposes the nonlinear transport problem into a sequence of linearized transport equations.
This section ends with the algorithm flowcharts of SLDG-RKEI with the adaptive time-stepping algorithm for the nonlinear Vlasov-Poisson and guiding center Vlasov systems.

\subsection{The SLDG scheme for linear transport problems}
\label{subsection:sldg}
We consider the general nonlinear transport equation in the form of \eqref{general_nonlinear} for  two-dimensional problems in a rectangular domain $\Omega$.
For the scope of this paper, we only consider the problems with the periodic boundary conditions.
In the special case when ${\bf P} (u;\mathbf{x},t)= (P_1(x, y, t), P_2(x, y, t))$ does not depend on $u$, the model problem \eqref{general_nonlinear} is linear and can be evolved by the SLDG scheme \cite{cai2016high} by accurately tracking characteristics with the velocity field ${\bf P}$.

\begin{figure}[h]
\centering
\subfigure[]{
\begin{tikzpicture}
    \draw[black,thin] (0,0.5) node[left] {} -- (5.5,0.5)
                                        node[right]{};
    \draw[black,thin] (0,2.) node[left] {$$} -- (5.5,2)
                                        node[right]{};
    \draw[black,thin] (0,3.5) node[left] {$$} -- (5.5,3.5)
                                        node[right]{};
    \draw[black,thin] (0,5 ) node[left] {$$} -- (5.5,5)
                                        node[right]{};
    \draw[black,thin] (0.5,0) node[left] {} -- (0.5,5.5)
                                        node[right]{};
    \draw[black,thin] (2,0) node[left] {$$} -- (2,5.5)
                                        node[right]{};
    \draw[black,thin] (3.5,0) node[left] {$$} -- (3.5,5.5)
                                        node[right]{};
    \draw[black,thin] (5,0) node[left] {$$} -- (5,5.5)
                                        node[right]{};
    \fill [blue] (3.5,3.5) circle (2pt) node[] {};
    \fill [blue] (5,3.5) circle (2pt) node[] {};
    \fill [blue] (3.5,5) circle (2pt) node[below right] {$A_j$} node[above left] {$c_4$};
    \fill [blue] (5,5) circle (2pt) node[] {};

     \draw[thick,blue] (3.5,3.5) node[left] {} -- (3.5,5)
                                        node[right]{};
      \draw[thick,blue] (3.5,3.5) node[left] {} -- (5,3.5)
                                        node[right]{};
       \draw[thick,blue] (3.5,5) node[left] {} -- (5,5)
                                        node[right]{};
        \draw[thick,blue] (5,3.5) node[left] {} -- (5,5)
                                        node[right]{};
    \fill [red] (1.,1) circle (2pt) node[above right,black] {};
    \fill [red] (3,1) circle (2pt) node[] {};
    \fill [red] (1,2.5) circle (2pt) node[below right] {$A_j^\star$} node[above left] {$c_4^\star$};
    \fill [red] (2.5,2.5) circle (2pt) node[] {};

     \draw[-latex,dashed](3.5,5)node[right,scale=1.0]{}
        to[out=240,in=70] (1,2.50) node[] {};

     \draw (0.5+0.01,2-0.01) node[fill=white,below right] {$A_l$};

     \draw [red,thick] (1,1)node[right,scale=1.0]{}
        to[out=20,in=150] (2,0.7) node[] {};

        \draw [red,thick] (2,0.7)node[right,scale=1.0]{}
        to[out=330,in=240] (3,1) node[] {};
             \draw [red,thick] (1,2.5)node[right,scale=1.0]{}
        to[out=310,in=90] (1.1,2) node[] {};
        \draw [red,thick] (1.1,2)node[right,scale=1.0]{}
        to[out=270,in=80] (1,1) node[] {};

        \draw [red,thick] (1,2.5)node[right,scale=1.0]{}
        to[out=10,in=180] (2.5,2.5) node[] {};

        \draw [red,thick] (3,1)node[right,scale=1.0]{}
        to[out=80,in=280] (2.5,2.5) node[] {};
\end{tikzpicture}
}
\subfigure[]{

\begin{tikzpicture}[scale = 1.3]
    \draw[black,thin] (0,0.5) node[left] {} -- (4,0.5)
                                        node[right]{};
    \draw[black,thin] (0,2.) node[left] {$$} -- (4,2)
                                        node[right]{};
    \draw[black,thin] (0,3.5) node[left] {$$} -- (4,3.5)
                                        node[right]{};
    \draw[black,thin] (0.5,0) node[left] {} -- (0.5,4)
                                        node[right]{};
    \draw[black,thin] (2,0) node[left] {$$} -- (2,4)
                                        node[right]{};
    \draw[black,thin] (3.5,0) node[left] {$$} -- (3.5,4)
                                        node[right]{};

    \fill [red] (1.,1) circle (2pt) node[above right,black] {$A_{j,l}^{\star}$};
    \fill [red] (3,1) circle (2pt) node[] {};
    \fill [red] (1,2.5) circle (2pt) node[below right] {$A_j^{\star}$} node[above left] {};
    \fill [red] (2.5,2.5) circle (2pt) node[] {};

     \draw (0.5+0.01,2-0.01) node[fill=white,below right] {};

\draw [red,thick] (1,1)node[right,scale=1.0]{}
        to[out=20,in=150] (2,0.7) node[] {};

        \draw [red,thick] (2,0.7)node[right,scale=1.0]{}
        to[out=330,in=240] (3,1) node[] {};
             \draw [red,thick] (1,2.5)node[right,scale=1.0]{}
        to[out=310,in=90] (1.1,2) node[] {};
        \draw [red,thick] (1.1,2)node[right,scale=1.0]{}
        to[out=270,in=80] (1,1) node[] {};

        \draw [red,thick] (1,2.5)node[right,scale=1.0]{}
        to[out=10,in=180] (2.5,2.5) node[] {};

        \draw [red,thick] (3,1)node[right,scale=1.0]{}
        to[out=80,in=280] (2.5,2.5) node[] {};
           \draw (0.5+0.01,2-0.01) node[fill=white,below right] {$A_l$};
         \draw[-latex,ultra thick] (1,1)node[right,scale=1.0]{}
        to  (2,1) node[] {};

             \draw[-latex,ultra thick]  (2,1)node[right,scale=1.0]{}
        to (2,2) node[] {};
             \draw[-latex,ultra thick]  (2,2)node[right,scale=1.0]{}
        to (1,2) node[] {};
             \draw[-latex,ultra thick]   (1,2)node[right,scale=1.0]{}
        to (1,1) node[] {};

\draw [thick] (1,1)-- (3,1) node[] {};
 \draw [thick] (1,2.5) -- (1,1) node[] {};
        \draw [thick] (1,2.5)--(2.5,2.5) node[] {};
        \draw [thick] (3,1)--(2.5,2.5) node[] {};
\end{tikzpicture}

}
\caption{Schematic illustration of the SLDG formulation in two dimension: quadrilateral approximation to an upstream cell.  }
\label{schematic_2d}
\end{figure}
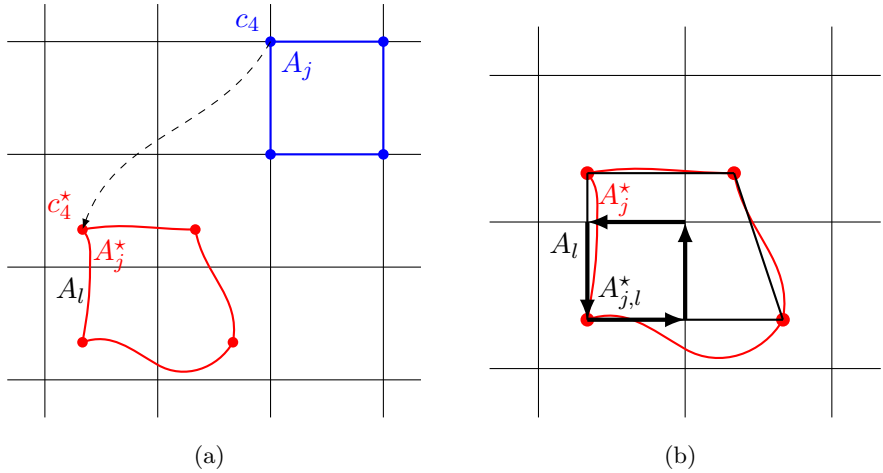

We partition the domain $\Omega$ by a set of non-overlapping tensor-product rectangular elements $A_j$, $j=1,\ldots,J$, and define the finite dimensional DG approximation space, $\bV_h^k = \{ v_h:  v_h|_{A_j} \in P^k(A_j) \}$, where $P^k(A_j)$ denotes the set of polynomials of degree at most $k$ over $A_j= [x_j^l,x_j^r]\times[y_j^b,y_j^t]$, where element centers and sizes are
$x_j = \frac{x_j^l + x_j^r}{2}$, $y_j = \frac{y_j^b + y_j^t }{2}$, $\Delta x_j = x_j^r - x_j^l$, $\Delta y_j = y_j^t - y_j^b$ respectively.
In the SLDG framework, we let the test function $\psi(x,y,t)$ satisfy the adjoint problem with   $\Psi\in P^k(A_j)$,
\begin{equation}
\begin{cases}
\psi_t + P_1(x,y,t) \psi_x + P_2(x,y,t) \psi_y = 0,\\
\psi(t=t_{n+1}) = \Psi.
\end{cases}
\end{equation}
Here we adopt the scaled Legendre polynomials.
For instance, for a $P^2$ polynomial,
 $\Psi$ varies the base in $\left\{1,\frac{x-x_j}{\Delta x_j}, \frac{y-y_j}{\Delta y_j},  \left( \frac{x-x_j}{\Delta x_j} \right)^2 -\frac{1}{12},
 \frac{ (x-x_j)(y-y_j) }{ \Delta x_j \Delta y_j },  \left( \frac{y - y_j }{\Delta y_j} \right)^2 -\frac{1}{12}  \right\}.$
Then, we have the identity
\begin{equation}
\frac{d}{dt} \int_{\widetilde{A}_{j}(t) } u(x,y,t) \psi(x,y,t) dxdy =0,
\label{2d_dt}
\end{equation}
where $\widetilde{A}_{j}(t)$ is the dynamic cell, moving from the Eulerian cell $A_{j}$ at $t^{n+1}$, i.e.,  $A_j=\widetilde{A}_{j}(t^{n+1})$, backward in time by following the characteristics trajectories, see eq.~\eqref{characteristic}. We denote $\widetilde{A}_{j}(t^n)$ as $A_{j}^\star$, i.e., the upstream cell bounded by the red curves in Figure \ref{schematic_2d}. The SLDG method is defined as follows. Given $u^n\in \bV_h^k$, we seek $u^{n+1}\in \bV_h^k$, such that for $\forall \Psi\in P^k(A_j)$, $j=1,\ldots,J$,
\begin{equation}
\int_{ A_{j}  } u^{n+1} \Psi dxdy =
\int_{ A_j^\star } u^n\psi(x,y,t^{n} ) dxdy.
\label{sl_2d}
\end{equation}
To update $u^{n+1}$, we need to properly evaluate the right-hand side (RHS) of \eqref{sl_2d}, the procedure of which we briefly review below.
In particular, we only review $P^1$ SLDG with quadrilateral approximation; to achieve third order accuracy, one can use the quadratic-curved quadrilateral approximation and $P^2$ polynomial solution space, see \cite{cai2016high} for more details regarding implementation.
For completeness, the procedure of the quadratic-curved approximation in given in Appendix \ref{append:b}.

\begin{enumerate}
\item
\noindent{\bf Characteristics tracing.}
Denote $c_q$, $q=1,\cdots,4$ as the four vertices of $A_j$  with the coordinates $(x_{j,q} , y_{j,q} )$.
We trace characteristics backward in time to $t^n$ for the four vertices
 by numerically solving the characteristics equations,
\begin{equation}
\begin{cases}
\frac{d x(t) }{dt} = P_1(x(t) ,y(t) ,t ),\\
\frac{ d y(t) }{dt} = P_2(x(t) ,y(t) ,t ), \\
x(t^{n+1} ) = x_{j,q},\\
y(t^{n+1} ) = y_{j,q},
\end{cases}
\label{characteristic}
\end{equation}
 and  obtain $c_q^\star$ with the new coordinate $( x_{j,q}^\star , y_{j,q}^\star ),\,q=1,\cdots,4$.
For example, see $c_4$ and $c_4^\star$ in Figure \ref{schematic_2d}.
In our implementation, a fifth order Runge-Kutta method \cite{butcher2008numerical} is used for solving \eqref{characteristic}.
 The upstream cell $A_j^\star$ can be approximated by a quadrilateral determined by the four vertices $c_q^\star$,
 which yields a second order approximation to sides of $A_j^\star$.

\item \noindent{\bf Evaluating the integrals over the upstream cells.}
Note that $u^n$ is a piecewise polynomial based on the partition. Then the integral over $A_j^\star$ has to be evaluated subregion-by-subregion. To this end, we denote $A_{j,l}^\star$ as a non-empty overlapping region between the upstream cell $A_j^\star$ and the background grid cell $A_l$, i.e., $A_{j,l}^\star = A_{j}^\star \cap A_l$, $A_{j,l}^\star \neq \emptyset$, and define the index set $\varepsilon_j^\star:=\{ l| A_{j,l}^\star \neq \emptyset \}$, see Figure \ref{schematic_2d} (b). The detailed procedure of detecting $A_{j,l}^\star$ can be found in \cite{cai2016high}.
The integral over the upstream cell $A_{j}^\star$ is broken up into the following integrals,
\begin{equation}
\int_{ A_{j}  } u^{n+1} \Psi dxdy
=
\sum_{l\in\varepsilon_j^\star } \int_{ A_{j,l}^\star } u^n\psi(x,y,t^{n} ) dxdy.
\label{temp1}
\end{equation}
Furthermore, $\psi(x,y,t^n)$ is not a polynomial in general, posing additional challenges for evaluating the integrals on the RHS of \eqref{temp1}. On the other hand, if the velocity field ${\bf P}$ is smooth, then $\psi(x,y,t^n)$ is smooth accordingly and can be well approximated by a polynomial. The following  procedure is then proposed.

 \begin{enumerate}
   \item  \textit{Least-squares approximation of test function $\psi(x,y,t^n)$.}
We use a least-squares procedure to approximate the test function $\psi(x,y,t^n)$  by a polynomial, based on the fact that  $\psi$ stays constant along characteristics. In particular, for $k=1$, we reconstruct a $P^1$ polynomial $\Psi^\star(x,y)$ by least-squares  with the interpolation constraints
\[
\Psi^\star(x_{j,q}^\star,y_{j,q}^\star) = \Psi(x_{j,q},y_{j,q}),\quad q=1,\ldots,4.
\]
Then,
\begin{equation}
\sum_{l\in\varepsilon_j^\star } \int_{ A_{j,l}^\star } u^n\psi(x,y,t^{n} ) dxdy \approx \sum_{l\in\varepsilon_j^\star } \int_{ A_{j,l}^\star } u^n\Psi^\star(x,y) dxdy.
\label{temp2}
\end{equation}


   \item  \textit{Line integral evaluation via Green's theorem.}
   Note that the integrands on the RHS of \eqref{temp2} are piecewise polynomials. To further simplify the implementation, we make use of  Green's theorem. We first introduce two auxiliary polynomial functions $P(x,y)$ and $Q(x,y)$ such that
   \begin{equation*}
   -\frac{\partial P }{\partial y } + \frac{\partial Q}{\partial x }  =  u(x,y,t^n)\Psi^\star(x,y).
   \end{equation*}
   Then area integral $ \int_{A_{j,l}^\star } u^n\Psi^\star(x,y)dxdy  $ can be converted into line integrals via Green's theorem, i.e.,
   \begin{equation}
   \int_{A_{j,l}^\star } u(x,y,t^n)\Psi^\star(x,y)dxdy = \oint_{\partial A_{j,l}^\star}  Pdx + Qdy,
   \label{Green}
   \end{equation}
   see Figure \ref{schematic_2d} (b). Note that the choices of $P$ and $Q$ are not unique, but the value of the line integrals is independent of the choices. These line integrals are organized into line integrals along inner segments and outer segments, for which numerical quadrature rules are applied to evaluate the line integrals, see Figure \ref{outer_inner}.

 \end{enumerate}

 \item
A positivity-preserving (PP) limiter  [52,7] is added to preserve the positivity of the solution, e.g. for the Vlasov-Poisson system. It can be implemented as follows.
 The numerical solution $u(x,y,t^n)$ in the cell $A_j$ is modified by $\tilde{u}(x,y)$,
 \begin{equation*}
 \tilde{u}(x,y) =
 \theta ( u(x,y,t^n) - \overline{u} ) + \overline{u},
 \theta = \min \{ \left| \frac{\overline{u}}{m'-\overline{u}}  \right|,1 \},
 \end{equation*}
 where $\overline{u}$ is the cell average of the numerical solution and $m'$ is the minimum value of $u(x,y,t^n)$ over $A_j$.
  Due to the PP limiter, the proposed conservative SLDG schemes can guarantee the $L^1$ stability of the Vlasov-Poisson solution, and the proof follows a similar argument in \cite{qiu2011positivity}.
\end{enumerate}

 \begin{rem}
 Note that the upstream cell may be traced to a location out of the computational domain.
For problems with the periodic boundary conditions, we define a set of ghost elements and fetch information from the other sides of the domain. For inflow type boundary condition, the SLDG scheme can be formulated accordingly over the boundary characteristics elements. General treatment of boundary conditions will be the subject of our future investigation.
 \end{rem}

 \begin{rem}
 To better approximate $\psi(x, y, t^n)$, one could use polynomials of degree higher than $k$ by sampling more traced back points. In this paper, we found the least square approximation by the polynomial with the same order $k$ sufficient.
 \end{rem}

\usetikzlibrary{shapes}

\begin{figure}[h!]
\centering
\subfigure[]{
\begin{tikzpicture}[scale = 1.]
    \draw[black,thin] (0,0.5) node[left] {} -- (4,0.5)
                                        node[right]{};
    \draw[black,thin] (0,2.) node[left] {$$} -- (4,2)
                                        node[right]{};
    \draw[black,thin] (0,3.5) node[left] {$$} -- (4,3.5)
                                        node[right]{};
    \draw[black,thin] (0.5,0) node[left] {} -- (0.5,4)
                                        node[right]{};
    \draw[black,thin] (2,0) node[left] {$$} -- (2,4)
                                        node[right]{};
    \draw[black,thin] (3.5,0) node[left] {$$} -- (3.5,4)
                                        node[right]{};

    \fill [red] (1.,1) circle (2pt) node[above right,black] {};
    \fill [red] (3,1) circle (2pt) node[] {};
    \fill [red] (1,2.5) circle (2pt) node[below right] {} node[above left] {};
    \fill [red] (2.5,2.5) circle (2pt) node[] {};

     \draw (0.5+0.01,2-0.01) node[fill=white,below right] {};
\draw [thick] (1,1)-- (3,1) node[] {};
 \draw [thick] (1,2.5) -- (1,1) node[] {};
        \draw [thick] (1,2.5)--(2.5,2.5) node[] {};
        \draw [thick] (3,1)--(2.5,2.5) node[] {};
   \usetikzlibrary{shapes.geometric}
  \node[fill,star,star points=4, star point ratio=.2] at (2,1) {};
  \node[fill,star,star points=4, star point ratio=.2] at (2,2.5) {};
  \node[fill,star,star points=4, star point ratio=.2] at (1,2) {};
  \node[fill,star,star points=4, star point ratio=.2] at (2.65,2) {};
  \draw [-latex] (1-0.2,2) -- node[right=3pt]{ }(1-0.2,1) node[] {};
  \draw [-latex] (1-0.2,2.5) -- (1-0.2,2) node[] {};

  \draw [-latex] (2.5,2.5+0.2) -- (2,2.5+0.2);
  \draw [-latex] (2,2.5+0.2) -- (1,2.5+0.2);

  \draw [-latex] (1,1-0.2) -- (2,1-0.2) node[] {};
  \draw [-latex] (2,1-0.2) -- (3,1-0.2) node[] {};
  \draw [-latex] (3+0.2,1) -- (2.65+0.2,2) node[] {};
  \draw [-latex] (2.65+0.2,2) -- (2.5+0.2,2.5) node[] {};



\end{tikzpicture}
}
\subfigure[]{
\begin{tikzpicture}[scale = 1.]
%

\node [below right,blue] at (2,1) {$s_1$};
\node [above,blue] at (2,2.5) {$s_2$};
\node [below,blue] at (0.9,2.) {$s_3$};
\node [above right,blue] at (2.6,2) {$s_4$};
\node [below right, blue] at (2,2) {$c_1$};
    \draw[black,thin] (0,0.5) node[left] {} -- (4,0.5)
                                        node[right]{};
    \draw[black,thin] (0,2.) node[left] {$$} -- (4,2)
                                        node[right]{};
    \draw[black,thin] (0,3.5) node[left] {$$} -- (4,3.5)
                                        node[right]{};
    \draw[black,thin] (0.5,0) node[left] {} -- (0.5,4)
                                        node[right]{};
    \draw[black,thin] (2,0) node[left] {$$} -- (2,4)
                                        node[right]{};
    \draw[black,thin] (3.5,0) node[left] {$$} -- (3.5,4)
                                        node[right]{};

    \fill [red] (1.,1) circle (2pt) node[above right,black] {};
    \fill [red] (3,1) circle (2pt) node[] {};
    \fill [red] (1,2.5) circle (2pt) node[below right] {} node[above left] {};
    \fill [red] (2.5,2.5) circle (2pt) node[] {};

   \usetikzlibrary{shapes.geometric}
  \node[fill,star,star points=4, star point ratio=.2,blue] at (2,1) {};
  \node[fill,star,star points=4, star point ratio=.2,blue] at (2,2.5) {};
  \node[fill,star,star points=4, star point ratio=.2,blue] at (1,2) {};
  \node[fill,star,star points=4, star point ratio=.2,blue] at (2.65,2) {};
  \node[fill,star,star points=4, star point ratio=.2,blue] at (2,2) {};

  \draw [-latex] (2-0.1,2+0.1) -- (2-0.1,2.5-0.1) node[] {};
  \draw [-latex] (2+0.1,2.5-0.1)--(2+0.1,2+0.1 )  node[] {};

   \draw [-latex] (2-0.1,1+0.1) --node[auto]{ } (2-0.1,2-0.1) node[] {};
  \draw [-latex] (2+0.1,2-0.1 )--(2+0.1,1+0.1 )  node[] {};
  \draw [-latex] (2.65-0.1,2-0.1) --(2+0.1,2-0.1)  node[] {};
  \draw [-latex] (2+0.1,2+0.1)--(2.65-0.1,2+0.1)  node[] {};
    \draw [-latex] (1+0.1,2+0.1) --(2-0.1,2+0.1)  node[] {};
  \draw [-latex] (2-0.1,2-0.1)--(1+0.1,2-0.1)  node[] {};

     \draw (0.5+0.01,2-0.01) node[fill=white,below right] {};
\draw [thick] (1,1)-- (3,1) node[] {};
 \draw [thick] (1,2.5) -- (1,1) node[] {};
        \draw [thick] (1,2.5)--(2.5,2.5) node[] {};
        \draw [thick] (3,1)--(2.5,2.5) node[] {};

\end{tikzpicture}
}

\caption{\footnotesize  Searching algorithm for outer (left) and inner (right) segments.}
\label{outer_inner}
\end{figure}
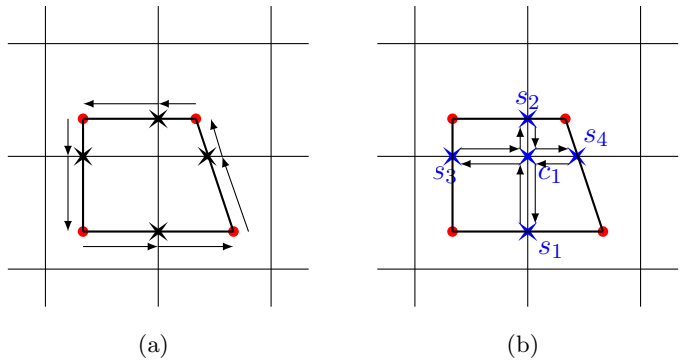

\subsection{SLDG for nonlinear models using Runge-Kutta exponential integrators}
For a general nonlinear model in the form of \eqref{general_nonlinear}, per time step evolution from $t^n$ to $t^{n+1}$, if we freeze the velocity field $\mathbf{P}(u)$ at $t^n$, the nonlinear problem is linearized around $u^n$ as follows
\begin{equation}
u_t + \nabla \cdot( \mathbf{P}( u^n ) u ) = 0.
\label{frozen_tn}
\end{equation}
One can apply a SL scheme (e.g. the SLDG method described above) to the linearized model \eqref{frozen_tn}, for which we adopt the notation of $SLDG(\mathbf{P}(u^n), \Delta t)$ for the update of solution from $u^n$ to $u^{n+1}$, i.e.
\begin{equation}
u^{n+1} = SLDG( \mathbf{P}(u^n), \Delta t) (u^n).
\label{eq: SL_1}
\end{equation}
Notice that due to the linearization of the velocity field $\mathbf{P}$ around $u^n$, there is a first order local truncation error in the temporal direction. To realize a high order temporal discretization, one could develop some strategies to track characteristics of the nonlinear dynamics, see \cite{cai2018high} for nonlinear Vlasov-Poisson and \cite{cai2019high} for the nonlinear guiding center Vlasov models. However, these strategies are problem-dependent and could be very algebraically involved for implementations. Here in this paper, we propose to adopt the RKEI scheme, which offers a unified framework and a black-box procedure for achieving high order temporal accuracy when coupled with the SLDG method for solving nonlinear problems.  We will first review the RKEI to solve the ordinary differential equations (ODEs) \cite{celledoni2003commutator}; then it is recognized that the exponential integrator for linear ODEs is equivalent to a semi-Lagrangian update of the solution for linear transport problems, if the spatial dimensions are kept continuous without numerical discretizations. With these observations, one can couple the high order RKEI with SLDG scheme for nonlinear dynamics.

\bigskip
\noindent
{\bf Review of the RKEI for nonlinear ODE systems.}
Consider a nonlinear ODE model of size $N$ in the form of
\begin{equation}
\frac{d Y(t) }{dt} = C(Y) Y, \quad Y(t=0) = Y_0,
\label{semi_nonlinear}
\end{equation}
where $C(Y)$ is a matrix-value function of size $N \times N$ that may depend on the solution $Y(t)$. In each time step, if we freeze $C(Y)$ at $C(Y^n)$, then we have a linearized problem $\frac{d Y(t) }{dt} = C(Y^n) Y \doteq C^n Y$, for which one can apply an exponential integrator to update the solution from $Y^n$ to $Y^{n+1}$ with a first order local truncation error:
\begin{equation}
Y^{n+1} = \exp(C^n \Delta t) Y^n.
\label{eq: exp_integrator_1}
\end{equation}
To improve the accuracy of the above first order scheme, a class of commutator-free exponential integrators can be used. The idea is to achieve high order temporal accuracy via taking composition of a sequence of linear solvers by freezing coefficients, which can be explicitly computed as a linear combination of $C(Y)$ from previous RK stages. In particular, the algorithm flowchart for the RKEI method is summarized as follows.

\bigskip
\begin{algorithm}[H]
  \caption{The commutator-free RKEI method \cite{celledoni2003commutator}.}
  \label{algo:inter}
  $p=Y^n$

  {\bf for} $r=1:s$ {\bf do}

   { \ \ \ \ \ \ } $Y_r = \text{exp}(\Delta t  \sum_k \alpha_{rJ^{(r)}}^k C(Y_k)  ) \cdots \text{exp}(\Delta t  \sum_k\alpha_{r1}^k C(Y_k)) p$

  {\bf end}

  $Y^{n+1} = \text{exp}( \Delta t \sum_k \beta_{J}^k  C(Y_k) ) \cdots \text{exp}(\Delta t  \sum_k\beta_1^k  C(Y_k) ) p $
  \end{algorithm}
Here $J^{(r)}$ represents the number of exponentials one has to take per RK stage.
The RKEI method can be represented by the Butcher tableau in Table \ref{table_rk}.
\begin{table}[htbp]
\centering
\begin{tabular}{r|c }
$\mathbf{c}$  &  $A$ \\ \hline
	          & $\mathbf{b}$\\
\end{tabular}
\caption{ A Butcher tableau for RKEI method, where  $a_{ik} = \sum_{l=1}^{J^{(i)}} \alpha_{i,l}^k$ and $b_k= \sum_{l=1}^{J^{n+1}} \beta_l^k$, which merges $J^{(i)}$  rows into one row in each stage $i$. }
\label{table_rk}
\end{table}
It is shown, in \cite{celledoni2009semi} by matching order conditions, that Butcher tableaus of many first and second order RK methods give the RKEI methods of the same order; but for third order RKEI method $J$ for some RK stages must be at least 2.
In this paper, we will also consider some other Butcher Tableaus from \cite{celledoni2003commutator, celledoni2009semi} as listed in Appendix \ref{append:a}.

\begin{exa}
A simple example is a first order exponential integrator, which is represented in Table \ref{CF1table}.
 \begin{table}[htbp]
\centering
\begin{tabular}{r|c } 
$0$&0  \\\hline
	&1\\
\end{tabular}
\caption{ CF1 \cite{pironneau1982transport}.  }
\label{CF1table}
\end{table}

\noindent It gives a simple first order method for the nonlinear model \eqref{general_nonlinear},
\begin{align}
Y^{(1)}  &  =  Y^n \nonumber \\
Y^{n+1}  &  =  \exp \left(   \Dt C(Y^{(1)})  \right) Y^n.
\end{align}
A third order RKEI method in \cite{celledoni2009semi} can be represented by the following Butcher tableau,
\begin{table}[htbp]
\centering
\begin{tabular}{r|c c c}
$0$             &                &               &  \\
$\frac12$    	& $\frac12$      &               &  \\
     1          & -1             &  2            & \\ \hline
                & $\frac{1}{12}$ & $\frac13$     & -$\frac14$\\
                & $\frac{1}{12}$  & $\frac13$      & $\frac{5}{12}$
\end{tabular}
\caption{ CF3G  }
\label{CF3Gtable}
\end{table}

\noindent
with which, the RKEI scheme for the nonlinear ODE system \eqref{semi_nonlinear} reads
\begin{align*}
Y^{(1)}  &  =  Y^n \\
Y^{(2)}  &  =  \exp \left( \frac12 \Dt C({Y^{(1)} })  \right) Y^n \\
Y^{(3)}  &  =  \exp \left( \Dt (-C({ Y^{(1)} }) + 2 C({ Y^{(2)} }))  \right) Y^n \\
Y^{n+1}  &  =  \exp \left( \Dt (\frac{1}{12} C({Y^{(1)} })  + \frac13 C({ Y^{(2)} }) + \frac{5}{12}  C({ Y^{(3)} })) \right) \\
 &\exp \left(\Dt( \frac{1}{12} C({ Y^{(1)} }) + \frac13  C({ Y^{(2)} }) - \frac14  C({  Y^{(3)}  }))   \right) Y^n.
\end{align*}
\end{exa}

\bigskip
\noindent
{\bf SLDG-RKEI for nonlinear transport problems.} It was recognized in \cite{celledoni2009semi} that the SL update of the solution \eqref{eq: SL_1} of a linearized transport problem is equivalent to applying the exponential integrator to the linearized ODE system \eqref{eq: exp_integrator_1}, as both of the methods evolve the differential equation {\em exactly} for a time step. Thus, the high order RKEI method developed in \cite{celledoni2003commutator} can be systematically coupled with the SLDG scheme described above for solving nonlinear dynamics \eqref{general_nonlinear} in a black-box manner.  In particular, a first order RKEI scheme with Butcher tableau \eqref{CF1table} gives rise to a first order scheme \eqref{eq: SL_1}; and a third order RKEI scheme with Butcher tableau \eqref{CF3Gtable} coupling with the SLDG scheme reads as
\begin{align}
u^{(1)}  &  =  u^n \nonumber \\
u^{(2)}  &  =  SLDG \left( \frac12 \mathbf{P}( u^{(1)} ),  \Dt    \right) u^n \nonumber \\
u^{(3)}  &  =  SLDG \left( - \mathbf{P}( u^{(1)} ) + 2 \mathbf{P}( u^{(2)} ) , \Dt    \right) u^n \nonumber \\
u^{n+1}  &  =  SLDG \left( \frac{1}{12} \mathbf{P}( u^{(1)} )  + \frac13 \mathbf{P}( u^{(2)} )    + \frac{5}{12} \mathbf{P}( u^{(3)} ),
\Delta t \right)
\nonumber \\
&SLDG \left( \frac{1}{12} \mathbf{P}( u^{(1)} )  + \frac13\mathbf{P}( u^{(2)} ) -   \frac14 \mathbf{P}( u^{(3)} ) , \Dt  \right) u^n.
\label{sl_cf3g}
\end{align}
The scheme is named as ``SLDG-CF3G" for short. The scheme is of high order accuracy in both space and time.

\subsection{Algorithm flowchart for nonlinear transport problems}
\label{section:adaptive}

Below we summarize the flowchart of the SLDG-RKEI method for the nonlinear transport problems, such as the nonlinear Vlasov-Poisson and the guiding center Vlasov systems. We borrow notations for RKEI schemes from Algorithm \ref{algo:inter}.

\bigskip
\begin{algorithm}[H]
  \caption{SLDG-RKEI algorithm to update solution from $t^n$ to $t^{n+1}$ for nonlinear transport problems in the form of \eqref{general_nonlinear}.}
  \label{algo:SLDG-RKEI}

  {\bf for} $r=1:s$ 

  { \ \ \ \ \ \ } $\bullet$  Let $u^{[1]}=u^n$, and $\mathbf{P}(u^{(r)})$ be the velocity field $\mathbf{P}(u)$  frozen at the stage $r$.

   { \ \ \ \ \ \ }   {\bf for} $l=1:J^{(r)}$



         { \ \ \ \ \ \ } { \ \ \ \ \ \ }$\bullet$
 $   u^{[l+1]} =SLDG  ( \sum_k \alpha_{rl}^k \textbf{P}(u^{(k)} ), \Delta t   ) u^{[l]}.  $

{ \ \ \ \ \ \ }   {\bf end}

 { \ \ \ \ \ \ } $\bullet$  Let $u^{(r)}=u^{[J^{(r)}+1]}$.

 {\bf end}

  $\bullet$  Let $u^{[1]}=u^n$.

   {\bf for} $l=1:J $

  { \ \ \ \ \ \ } $\bullet$
 $   u^{[l+1]} =SLDG  ( \sum_k \beta_{l}^k \textbf{P}(u^{(k)} ), \Delta t   ) u^{[l]}.  $

    {\bf end}

  $\bullet$ $u^{n+1}=u^{[J+1]}$.

  \end{algorithm}
We make the following two remarks on the flow chart for the Vlasov-Poisson system \eqref{vlasov} and the guiding center Vlasov model \eqref{guiding} applications.

1. To freeze the velocity field, one can apply an LDG method in \cite{arnold2002unified, cockburn1998local, castillo2000priori,rossmanith2011positivity} for \eqref{poisson} and \eqref{poisson2d}, respectively. The resulting velocity fields are discontinuous. In this paper, we take the average of the velocity field in tracking characteristics. A post-processing technique \cite{ryan2005extension} could be added for smoothing the discontinuous velocity field, for which we plan to explore in the future.


%
%
%
%
%
%


2.
In \cite{cai2019high}, an adaptive time-stepping algorithm is introduced for the systems \eqref{vlasov} and \eqref{guiding}, based on controlling the $L^\infty$ norm of relative deviation of upstream cells' area. This adaptive time-stepping algorithm is summarized as Algorithm \ref{algo:adaptive} below and can be applied to each RK stage of the SLDG-RKEI algorithm in order to adaptively choose time stepping sizes.

\bigskip
\begin{algorithm}[H]
  \caption{Adaptive Time-Stepping Algorithm in an SLDG-RKEI method.}
  \label{algo:adaptive}
  \begin{description}
\item
  Let $A_j$ be an Eulerian cell and $A_j^\star$ be its upstream cell in a stage of the SLDG-RKEI scheme.
\item
  Compute $\theta= \max_j \left| \frac{ \text{area}\left(A_j^{\star}\right) - \text{area}\left(A_j\right)  }{ \text{area}\left(A_j\right)} \right|$.
  \item
 Let $\delta_M$ and $\delta_m$ be prescribed thresholds for decreasing and increasing $CFL$ number. In our simulations,  $\delta_M=1\%$ and $\delta_m=0.05\%$.
       \begin{description}
         \item[if] $\theta>\delta_M$, {\bf then} $irefine = 1$
            we decrease $CFL$ (e.g. $CFL= \max(CFL-1,CFL_{\min})$),  and restart the SLDG-RKEI scheme with the updated CFL.
         \item[else if] $\theta<\delta_m$, {\bf then} $irefine =0$,
            we increase $CFL$ (e.g. $CFL=\min( CFL+1, CFL_{\max})$) and restart the SLDG-RKEI scheme with the updated CFL.
         \item[else] Continue.
         \item[end if]
       \end{description}
       \end{description}
  \end{algorithm}

\section{Numerical results}
\label{sec3}

\subsection{Nonlinear Vlasov dynamics}
In this section, we present numerical results of the SLDG-RKEI for the nonlinear Vlasov-Poisson system and the Guiding center Vlasov model.
Unless otherwise noted, we use the following notations: the SLDG method with $P^k$ polynomial basis is denoted as $P^k$ SLDG;
we use the notation without or with -QC to denote  quadrilateral or quadratic-curved (QC) quadrilateral approximation to upstream cells.
We set the time step as
\begin{equation}
\Delta t = \frac{CFL}{ \frac{a}{\Delta x}+\frac{b}{\Delta y} } ,
\end{equation}
where $a$ and $b$ are maximum transport speeds in $x$- and $y$- directions, respectively.
 For  CPU  time  comparison,  all  simulations  are  performed  on  a computer with  Intel(R) Xeon(R) CPU E5-2660 v3 @ 2.60GHz.

\begin{exa}
\emph{(VP system: strong Landau damping.)} Consider strong Landau damping for the VP system with
the initial condition being a perturbed equilibrium
\begin{equation}
f(x,v,t=0) = \frac{ 1 }{\sqrt{2\pi} } ( 1+ \alpha \cos(kx)  ) \exp \left( -\frac{v^2}{2} \right) ,
\label{init}
\end{equation}
where $\alpha = 0.5$ and $k=0.5$.
The computational domain is $[0,4\pi]\times[-2\pi,2\pi]$. This problem has been numerically investigated by several authors, e.g. see \cite{xiong2014high,zhu2016h,huang2016semi}.

In Table \ref{VP_spatial}, we first test the spatial convergence of the SLDG methods with the third order temporal scheme, CF3C03, whose Butcher tableau is presented in the Appendix \ref{append:a}.
The well-known time reversibility of the VP system \cite{degond2004modeling} is used to test the order of convergence.
In particular, one can integrate the VP system forward to some time $T$, and then reverse
the velocity field of the solution and continue to integrate the system by the same amount of time $T$. Then, the solution should recover the initial condition with reverse velocity field, which can be used as a reference solution.
We show the $L^1$, $L^2$, and $L^\infty$ errors and the corresponding orders of convergence for $P^k$ SLDG(-QC)+CF3C03 scheme, $k= 1,2$ with $CFL=0.1$ in Table \ref{VP_spatial}.
We observe the second order convergence  for $P^1$ SLDG scheme; we observe a second order convergence  for $P^2$ SLDG scheme with quadrilateral approximation to upstream cells, and a third order convergence for the $P^2$ SLDG-QC scheme.

\begin{table}[!ht]
\caption{ Strong Landau damping.  $T=0.5$. Use the time reversibility of the VP system. Order of accuracy in space for $P^k$ SLDG(-QC)+CF3C03 scheme, $k= 1,2$. We set $CFL=0.1$ so that the spatial error is the dominant error.
}
\vspace{0.1in}
\centering
\begin{tabular}{ l  cc cc cc }
\hline
Mesh  &{$L^1$ error} & Order    &{$L^2$ error} & Order   &{$L^\infty$ error} & Order  \\
 \hline
  \multicolumn{7}{l}{ $P^1$ SLDG+CF3C03}
     \\
  $32^2$ &     5.88E-04 & &     1.21E-03 & &     1.18E-02 & \\
  $64^2$ &     1.50E-04 &    1.97 &     3.17E-04 &    1.94 &     3.49E-03 &    1.76 \\
  $96^2$ &     6.67E-05 &    1.99 &     1.42E-04 &    1.97 &     1.61E-03 &    1.91 \\
  $128^2$&     3.76E-05 &    2.00 &     8.06E-05 &    1.98 &     9.19E-04 &    1.95 \\
  $160^2$ &     2.41E-05 &    2.00 &     5.17E-05 &    1.99 &     5.93E-04 &    1.97\\

   \multicolumn{7}{l}{ $P^2$ SLDG+CF3C03}
     \\
   $32^2$ &     9.59E-05 & &     2.18E-04 & &     1.95E-03 & \\
   $64^2$ &     2.43E-05 &    1.98 &     5.57E-05 &    1.97 &     5.05E-04 &    1.95  \\
   $96^2$ &     1.09E-05 &    1.98 &     2.50E-05 &    1.98 &     2.26E-04 &    1.99 \\
   $128^2$ &    6.15E-06 &    1.99 &     1.41E-05 &    1.98 &     1.27E-04 &    1.99 \\
   $160^2$ &     3.94E-06 &    1.99 &     9.07E-06 &    1.99 &     8.15E-05 &    2.00 \\

   \multicolumn{7}{l}{ $P^2$ SLDG-QC+CF3C03}
     \\
    $32^2$ &     3.69E-05 & &     8.39E-05 & &     1.09E-03 & \\
    $64^2$ &     4.39E-06 &    3.07 &     1.03E-05 &    3.03 &     1.38E-04 &    2.98 \\
   $96^2$ &     1.28E-06 &    3.04 &     3.02E-06 &    3.02 &     4.09E-05 &    3.00 \\
  $128^2$ &     5.37E-07 &    3.02 &     1.27E-06 &    3.01 &     1.71E-05 &    3.03 \\
   $160^2$ &     2.74E-07 &    3.02 &     6.50E-07 &    3.01 &     9.97E-06 &    2.42 \\
\hline

\end{tabular}
\label{VP_spatial}
\end{table}

We then test the temporal convergence of different  temporal schemes by the strong Landau damping test case integrated to $T=5$.
In order to minimize the errors from the spatial discretization, we adopt the $P^2$ SLDG-QC scheme with a fixed mesh of $160\times160$ cells.
The reference solution is computed by the $P^2$ SLDG-QC scheme with the same mesh but using a small $CFL=0.1$.
We report plots of $L^1$  error  versus the $CFL$ number and $L^1$ error versus CPU cost of different RKEI methods in Figure \ref{VP_time}.
We make the following observations:
(1) Expected temporal orders of convergence are observed for various SLDG-RKEI method.
There is a plateau region for third order temporal schemes for CFL ranging from 5 to 10, for which we think is due to the interaction between spatial and temporal errors.
(2) $CFL$s can be taken to be as large as $50$, which is much larger than that for an Eulerian RKDG method, whose $CFL$ upper bound is $1/(2k+1)$ with $k$ being the degree of the polynomial.
(3) The error versus CPU plot is closely related to the error versus CFL plot, as the CFL number is linearly related to the number of time steps taken. When the CFL number is relatively large and the  temporal errors dominate, we find that the schemes with higher-order temporal accuracy more efficient than the lower order ones;
 (4) For comparison, we plot results from the SLDG method coupled with a predictor-corrector way of tracking characteristics \cite{cai2018high}. It is observed that third order temporal accuracy is numerically achieved for all third order time integrators; error magnitudes from SLDG-RKEI schemes are observed to be smaller.
 (5) To show the high order of convergence of the SLDG without operator splitting, the results of the $Q^2$ SLDG method with the Strang splitting method in \cite{cai2019comparison} are presented.  We observe that $Q^2$ SLDG-split is subject to the second order splitting error, while we observe the third order of convergence for the proposed schemes with the third order temporal methods. On the other hand, the CPU cost of the splitting algorithm is lower than that of the 2D algorithm.
 
\begin{figure}[h]
\centering                              
\includegraphics[height=70mm]{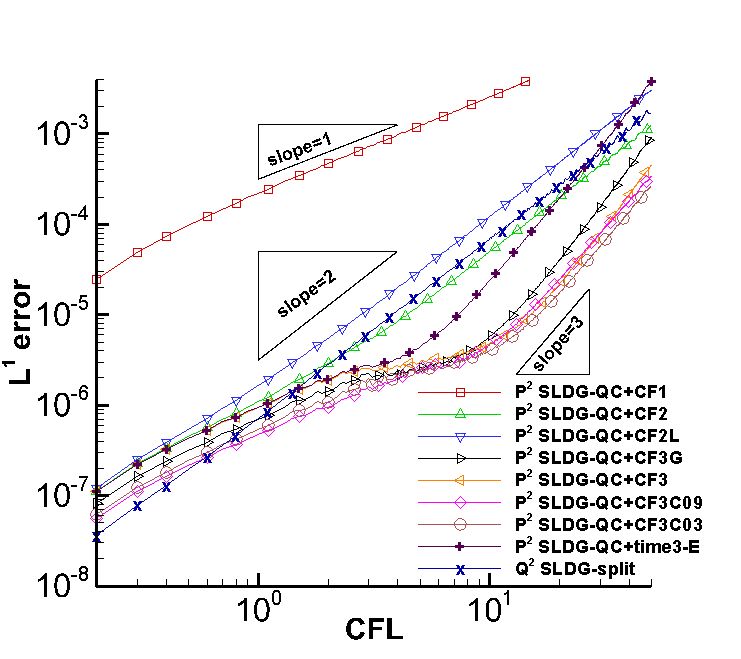}
\includegraphics[height=70mm]{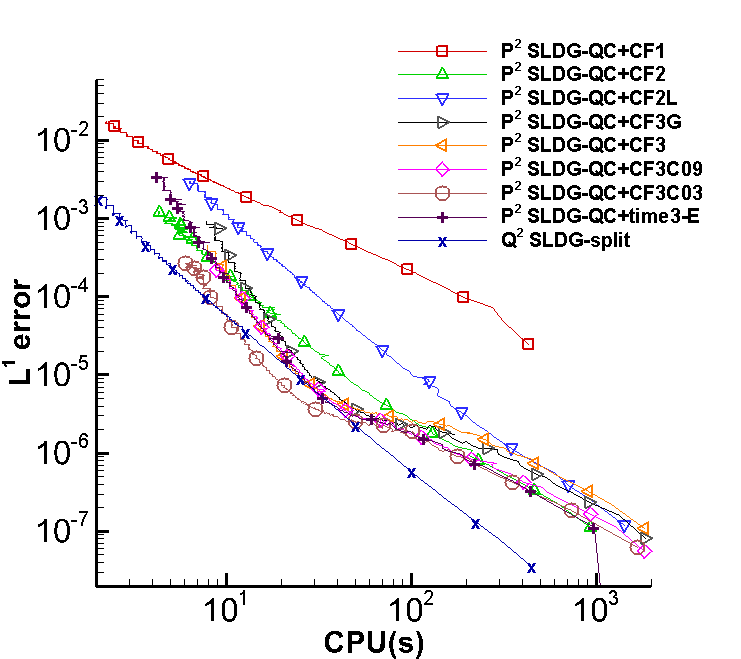}
\caption{ Plots of error versus the $CFL$ number (left) and error versus CPU (s) (right) for solving strong Landau damping at $T=5$.
Temporal order of convergence in $L^1$ norm of $P^2$ SLDG-QC scheme with various temporal schemes (denoted as $P^2$ SLDG-QC+temporal scheme) and $Q^2$ SLDG-split by comparing numerical solutions with
a reference solution from the corresponding scheme with $CFL = 0.1$.   The mesh of $160\times160$ is used.
$P^2$ SLDG-QC+time3-E is $P^2$ SLDG-QC with the third order prediction correction method in \cite{cai2018high}. $Q^2$ SLDG-split is the $Q^2$ SLDG with the Strang splitting in \cite{cai2019comparison}.
The Butcher tableaus of CF1, CF2, CF2L, CF3G, CF3, CF0C09, CF3C03 can be found in the previous section and the Appendix \ref{append:a}.   }
\label{VP_time}
\end{figure}

There are several conserved quantities in the VP system which should remain constant in time. These include the $L^p$ norm, kinetic energy and entropy:
\begin{itemize}
  \item $L^p$ norm, $1\leq p \leq \infty$:
      \begin{equation}
        \| f \|_p = \left(  \int_v \int_x |f(x,v,t)|^p\,dxdv \right)^{\frac{1}{p}},
      \end{equation}
  \item Energy:
      \begin{equation}
        \text{Energy} = \int_v \int_x f(x,v,t) v^2 dx dv + \int_x E^2 (x,t)dx,
      \end{equation}
  \item Entropy:
     \begin{equation}
      \text{Entropy} = \int_v \int_x f(x,v,t) \log (f(x,v,t) ) dxdv.
     \end{equation}
\end{itemize}
In Figure \ref{VP_norms}, we plot the time evolution of the relative deviation of $L^1$ and $L^2$  norms of the solution as well as the discrete kinetic energy and entropy.
Here we choose to present only a few representative RKEI schemes and run our simulations with $CFL=10$.
A few observations can be made:
(1) due to the truncation of the velocity domain, the error for the relative deviation of $L^1$ norm is on the of $10^{-9}$;
(2) the SLDG with higher degree polynomial does a better job in conserving these physical norms than the SLDG with lower degree polynomial;
(3) there is little difference in the performance of preserving norms for the $P^k$ SLDG with various same order temporal scheme including the prediction-correction method.
In Figure \ref{VP_surface}, we present the surface plots of the solutions at $T=40$ computed by $P^k$ SLDG ($k=1,2$) with the mesh of $160\times160$ elements.
The numerical solution of $P^2$ SLDG outperform that of $P^1$ SLDG in terms of resolution.

\begin{figure}[h]
\centering                              
\includegraphics[height=60mm]{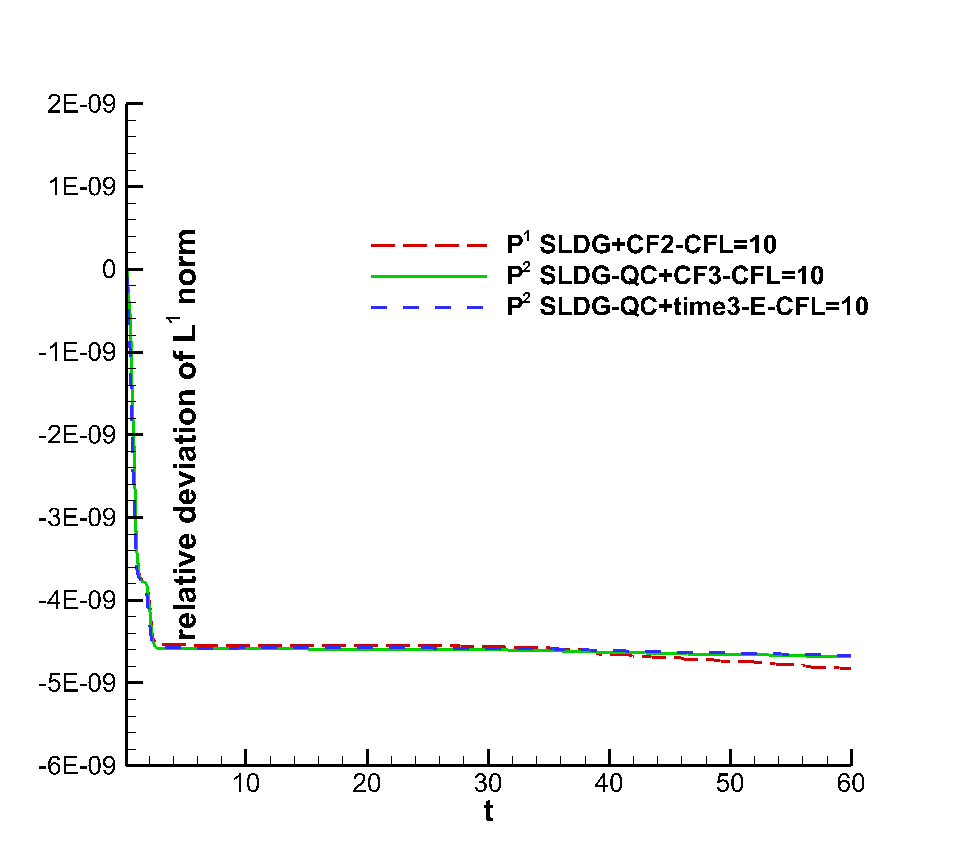}
\includegraphics[height=60mm]{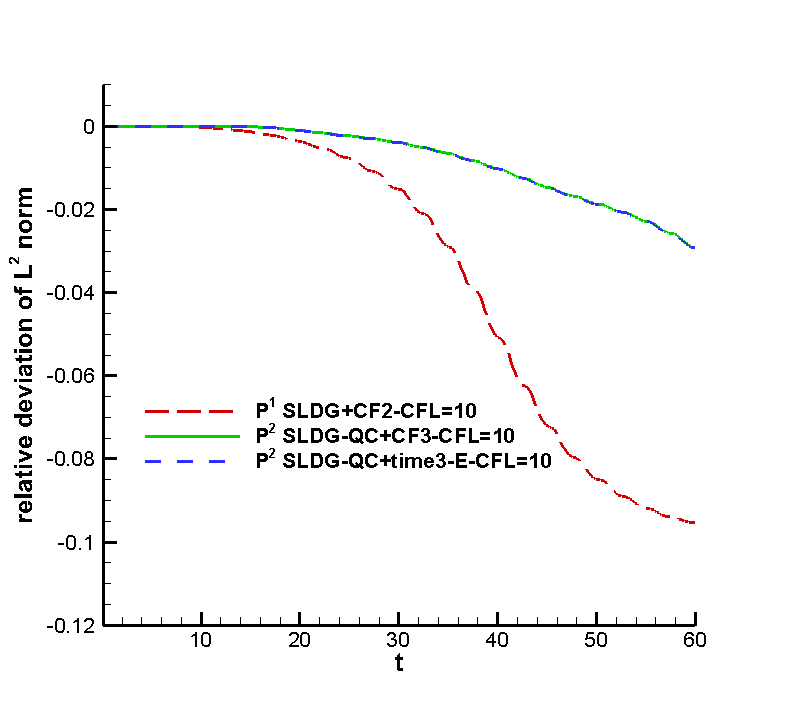}
\includegraphics[height=60mm]{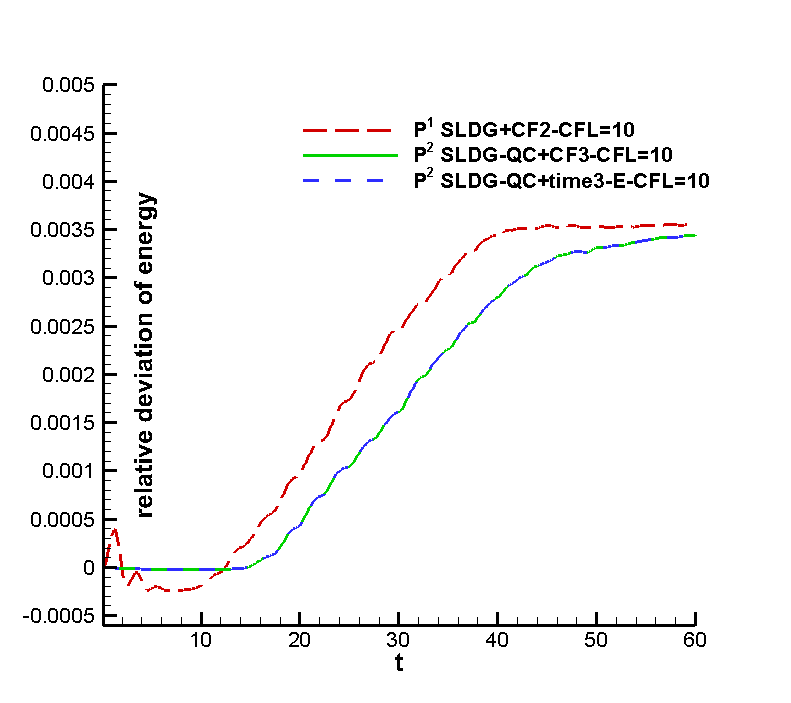}
\includegraphics[height=60mm]{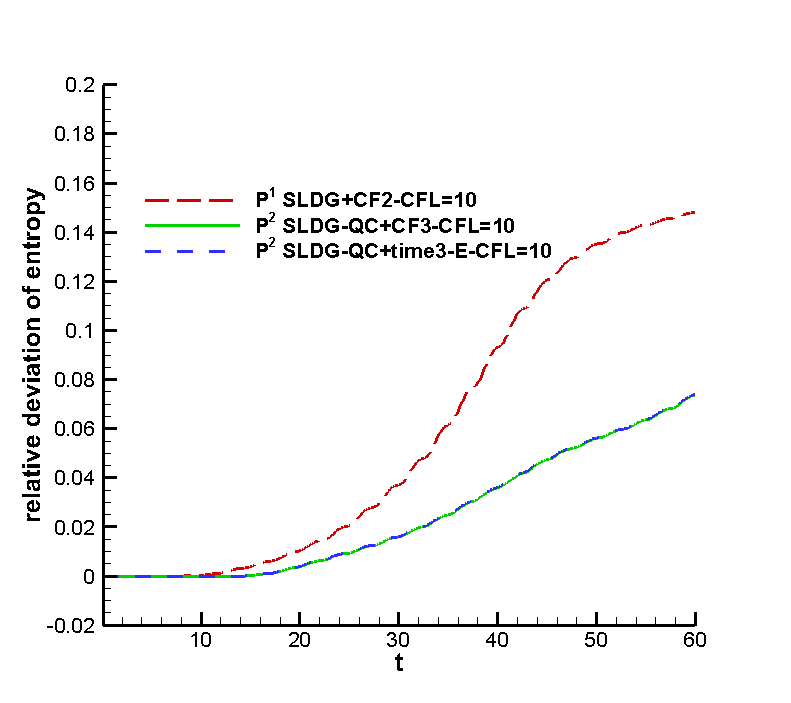}
\caption{Strong Landau damping. Time evolution of the relative deviation of $L^1$ (upper left) and $L^2$ (upper right) norms of the solution as well as the discrete kinetic energy (lower left) and entropy (lower right). We use a mesh of $160\times160$ cells and  $CFL=10$. $P^2$ SLDG-QC+time3-E is $P^2$ SLDG-QC with the third order prediction correction method in \cite{cai2018high}. }
\label{VP_norms}
\end{figure}

\begin{figure}[h]
\centering                              
\includegraphics[height=60mm]{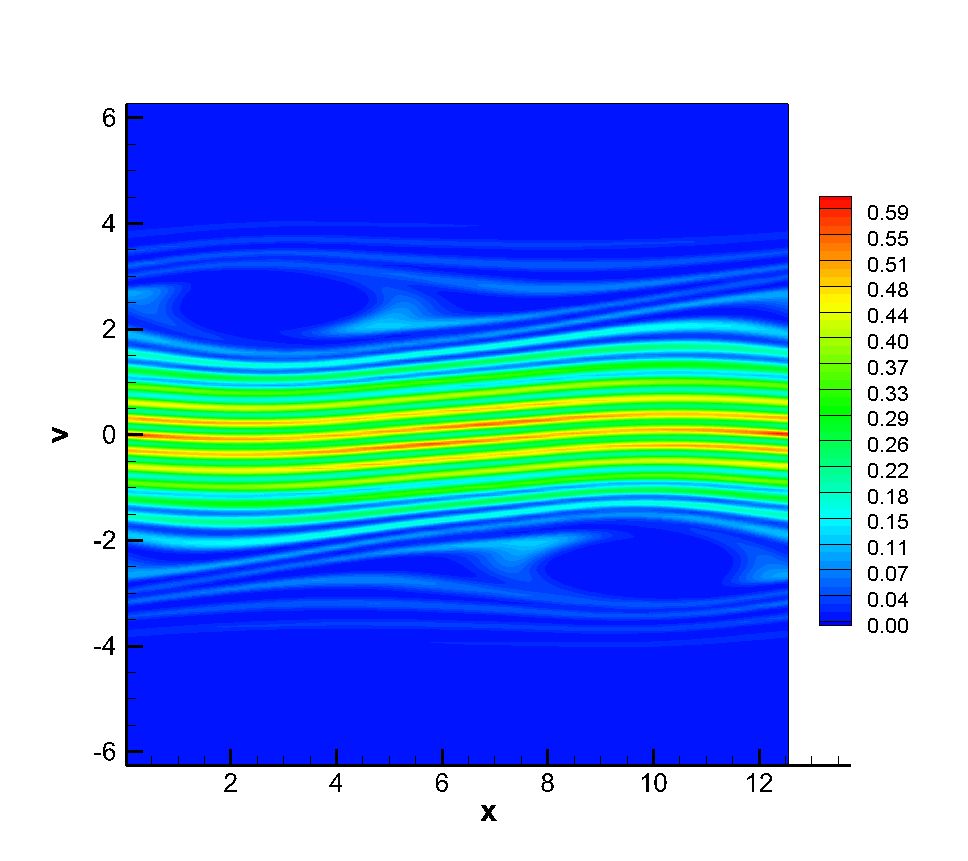}
\includegraphics[height=60mm]{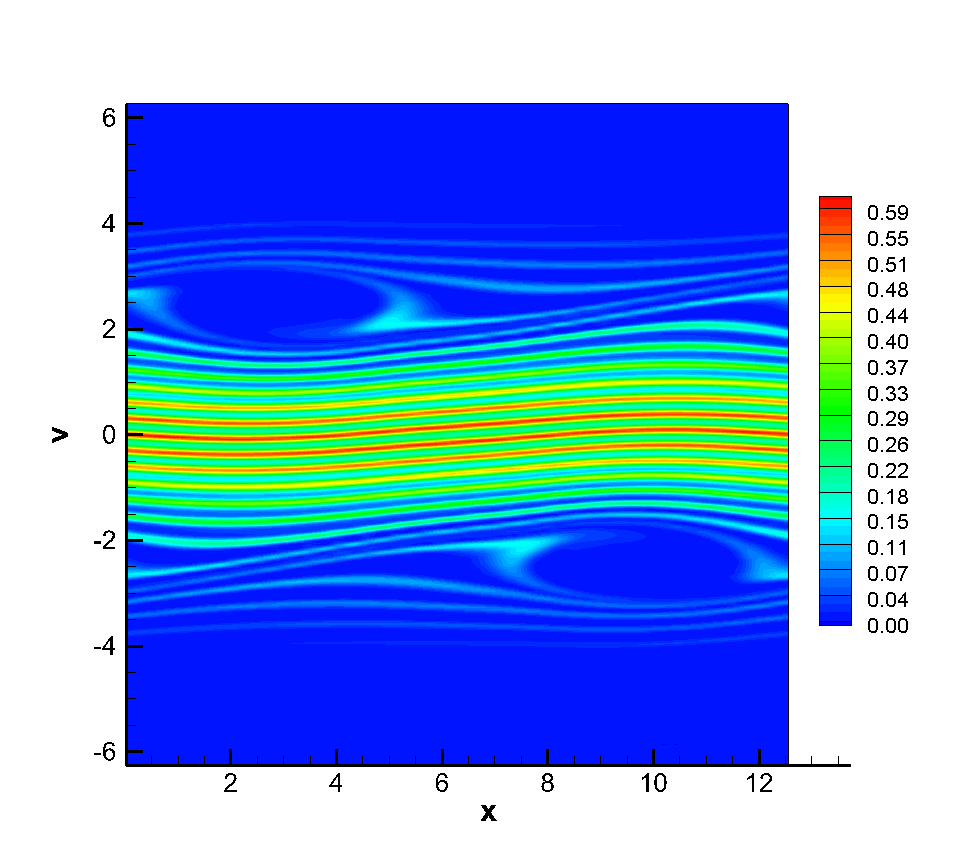}

\caption{Surface plots of the numerical solutions for the strong Landau damping at $T=40$.
We use a mesh of $160\times160$ cells and $CFL=10$.
left: $P^1$ SLDG+CF2;
right: $P^2$ SLDG-QC+CF3C03. }
\label{VP_surface}
\end{figure}

\end{exa}

\begin{exa}
 
\emph{(VP system:  Landau damping.)}
In this example, we assess the optimal convergence and superconvergence property of the 2D SLDG scheme. We first consider weak Landau damping for the VP system. The initial condition is the same as the strong one \eqref{init}, but with a smaller perturbation parameter $\alpha=0.01$. In this example, we show the optimal and superconvergence of SLDG scheme for the VP system. We first recall the optimal convergent and superconvergent rates of SLDG for the 1D linear transport problem in \cite{yang2020optimal}. It is shown that, if the solution is smooth enough, the $L^2$ norm of the SLDG error has the following estimate,
\begin{equation}
\label{eq: supe}
 \|u-u_h\|_{L^2} \le
 C_1 h^k+C_2 h^{2k+1}t,
 \end{equation}
where $h$ and $k$ are the mesh size and the approximation order for spatial discretization, $C_1$ and $C_2$ are constants independent of the mesh size and $t$ is the integration time. The first term in \eqref{eq: supe} is the projection error that is time independent and is the dominant error until time $t$ becomes very large, while the second term is time dependent and has a super convergent rate of $2k+1$. Such phenomenon can be observed for the weak Landau damping test by repeatedly first going forward in time with time stepping size determined by $CFL=0.5$ for a given mesh to $T=0.2$ and then backward in time returning to $T=0$, due to time reversibility of the VP system \cite{degond2004modeling}.
We adopt such strategy in the test so that the solution remains smooth enough, which is the key assumption of the error estimate. We test $P^1$ SLDG scheme with a third order temporal scheme CF3C03, with uniform meshes $N=20^2,30^2,40^2$. In Figure \ref{longtime}, we plot the errors of the solutions of the $P^1$ SLDG scheme at time $t=0.4n, n=1,2,3,\cdots$.
 It is observed that (1) the error of the SLDG solution does not grow in time after very long time; (2) when the error starts to grow, we report the slope of the error with respect to time (i.e. $C_2 h^{2k+1}$ according to \eqref{eq: supe}) in the Table on the right-hand panel of Figure \ref{longtime}.  These slopes are measured by the divided difference of error and time around the region where the error shows linear growth with time. Specifically, $\text{slope}= \frac{error(t_2)-error(t_1)}{t_2 -t_1} $, where $(t_1,t_2)=(2000,3000),(t_1,t_2)=(3000,5000),(t_1,t_2)=(4000,8000)$ for estimating the slopes for $N=20^2,30^2,40^2$ respectively.
 Third order is observed from these slopes and is consistent with the error estimate in \eqref{eq: supe}. On the other hand, we would like to point out that the superconvergence result is observed for VP system only at the ``linear'' regime.  When the nonlinear effect becomes significant, e.g. the strong Landau damping, the superconvergent property is lost.
 
  \begin{figure}[h!]
    \centering
    \includegraphics[height=60mm]{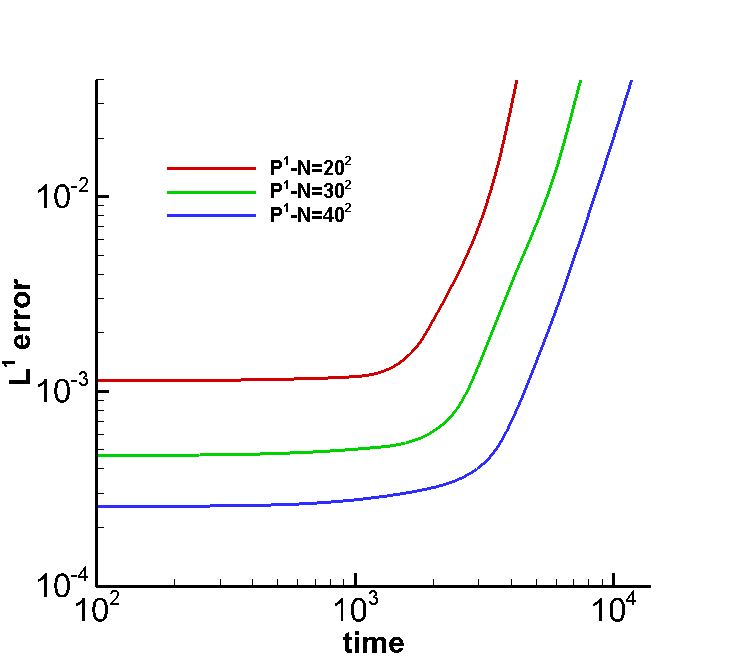}
    \qquad
\begin{tabular}[b]{ cc c  }
\hline
{ N  }  &{Slope } & {Order}      \\
\hline
\multicolumn{3}{c}{ $P^1$}\\ \hline
    $20^2$ &     $1.32\times10^{-5}$ &        \\
    $30^2$ &    $4.23\times10^{-6}$ &  2.81    \\
    $40^2$ &    $1.88\times10^{-6}$ & 2.83    \\
    \hline
 \\
 \multicolumn{3}{l}{ \quad}
\end{tabular}
    \caption{
    Weak Landau damping.
    The solution is evolved by repeatly evolving the solution first going forward in time with $T = 0.2$ and then backward in time returning to $T = 0$. The errors are tracked at $t=0.4n, n=1,2,3,\cdots$. $P^1$ SLDG+CF3C03 with $CFL=0.5$.
      Uniform  meshes have $N=20^2, 30^2, 40^2$ elements.   }
    \label{longtime}
  \end{figure}
\end{exa}
\begin{exa}
(The guiding center Vlasov system: spatial accuracy and convergence test).
Consider the guiding center Vlasov model on the domain $[0,2\pi]\times[0,2\pi]$ with the initial condition
\begin{equation}
\rho (x,y,0) = -2 \sin(x) \sin(y)
\end{equation}
and the periodic boundary condition.
The exact solution stays stationary.
We test the spatial convergence of the proposed SLDG methods with the third order temporal scheme, CF3C03, for solving  the guiding center Vlasov model up to time $T=1$ and report these results in Table \ref{Euler_spatial}.
We observe the expected second order of convergence for $P^1$ SLDG+$P^2$ LDG in $L^1$, $L^2$, and $L^\infty$ norms.
For $P^2$ SLDG scheme with quadrilateral approximation to upstream cells, the error magnitude will be reduced but still  second order.
For $P^2$ SLDG scheme with quadratic-curved quadrilateral approximation to upstream cells, we observe the third order of convergence in $L^1$ and $L^2$ norm but the second order of convergence in $L^\infty$ norm.


\begin{table}[!ht]
\caption{
The guiding center Vlasov system on the domain $[0,2\pi]\times[0,2\pi]$ with the initial condition
$
\omega (x,y,0) = -2 \sin(x) \sin(y).
$
Periodic boundary conditions in two directions.
$T=1$.
$CFL=1$.
 }
\centering
\begin{tabular}{ccccc cc  }
\hline
{  Mesh }  &{$L^1$ error} & Order    &{$L^2$ error} & Order  & {$L^\infty$ error} & Order  \\
\hline
%

  \multicolumn{7}{l}{ $P^1$ SLDG+$P^2$ LDG+CF3C03  }
\\
    $20^2$ &     1.39E-02 & &     1.88E-02 & &     1.06E-01 & \\
    $40^2$ &     3.66E-03 &     1.93 &     4.97E-03 &     1.92 &     3.12E-02 &     1.76 \\
    $60^2$ &     1.65E-03 &     1.97 &     2.24E-03 &     1.97 &     1.44E-02 &     1.90 \\
    $80^2$ &     9.37E-04 &     1.96 &     1.27E-03 &     1.95 &     8.27E-03 &     1.93 \\
   $100^2$ &     6.01E-04 &     1.99 &     8.17E-04 &     1.99 &     5.34E-03 &     1.96 \\

   \multicolumn{7}{l}{ $P^2$ SLDG+$P^3$ LDG+CF3C03 }
\\

    $20^2$ &     4.52E-03 & &     6.61E-03 & &     5.60E-02 & \\
    $40^2$ &     1.02E-03 &     2.14 &     1.53E-03 &     2.11 &     1.49E-02 &     1.91 \\
    $60^2$ &     4.30E-04 &     2.14 &     6.37E-04 &     2.17 &     6.76E-03 &     1.95 \\
    $80^2$ &     2.54E-04 &     1.82 &     3.80E-04 &     1.79 &     3.89E-03 &     1.92 \\
   $100^2$ &     1.52E-04 &     2.30 &     2.29E-04 &     2.27 &     2.52E-03 &     1.94 \\

   \multicolumn{7}{l}{ $P^2$ SLDG-QC+$P^3$ LDG+CF3C03 }
\\

    $20^2$ &     2.13E-03 & &     2.77E-03 & &     2.06E-02 & \\
    $40^2$ &     2.73E-04 &     2.97 &     3.63E-04 &     2.93 &     4.72E-03 &     2.13 \\
    $60^2$ &     8.11E-05 &     2.99 &     1.09E-04 &     2.96 &     2.06E-03 &     2.04 \\
    $80^2$ &     3.48E-05 &     2.94 &     4.74E-05 &     2.91 &     1.14E-03 &     2.05 \\
   $100^2$ &     1.77E-05 &     3.02 &     2.44E-05 &     2.98 &     7.28E-04 &     2.02 \\
   \hline

\end{tabular}
\label{Euler_spatial}
\end{table}

\end{exa}


\begin{exa}
(The guiding center Vlasov model: Kelvin-Helmholtz instability problem).
This is the two-dimensional guiding center model problem \eqref{guiding} with the initial condition
\begin{equation}
\rho_0(x,y) = \sin(y)+0.015 \cos(kx)
\end{equation}
and periodic boundary conditions on the domain $[0,4\pi]\times[0,2\pi]$.
We let $k=0.5$, which will create a Kelvin-Helmholtz instability \cite{shoucri1981two}, which is well studied numerically by many authors before (e.g. see \cite{zhu2017h,cai2019high}).

First, we  test the temporal convergence  of the proposed SLDG method with different temporal schemes by computing this problem up to $T=5$.
In order to minimize the errors from the spatial scheme, we adopt the $P^2$ SLDG-QC scheme with $P^3$ LDG method using a fixed mesh of $120\times120$ cells.
The reference solution is computed by the same scheme with the same mesh but using a small $CFL=0.1$.
 We report plots of  $L^1$ error  versus the CFL number and $L^1$ error  versus CPU cost in Figure \ref{KH_time}. 
Some observations can be concluded from Figure \ref{KH_time}:
(1) expected order of convergence is observed for all temporal schemes; and $CFL$ number can be taken to be as large as $50$;
(2) by comparing the error magnitudes, CF2L performs better than CF2 and CF3C09 performs the best;
 (3) by comparing CPU cost, we find that the $P^2$ SLDG-QC schemes the third-order temporal accuracy  much more efficient  than the lower-order ones.

\begin{figure}[h]
\centering                              
\includegraphics[height=70mm]{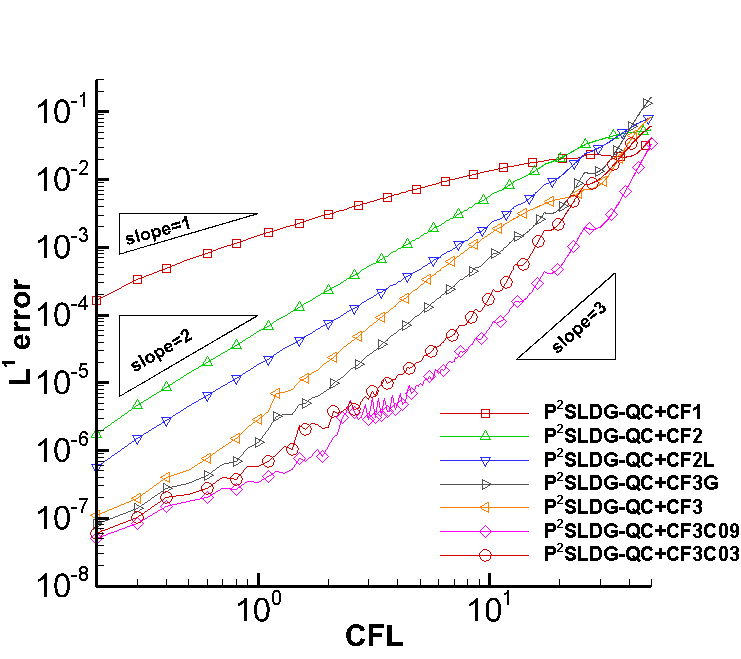}
\includegraphics[height=70mm]{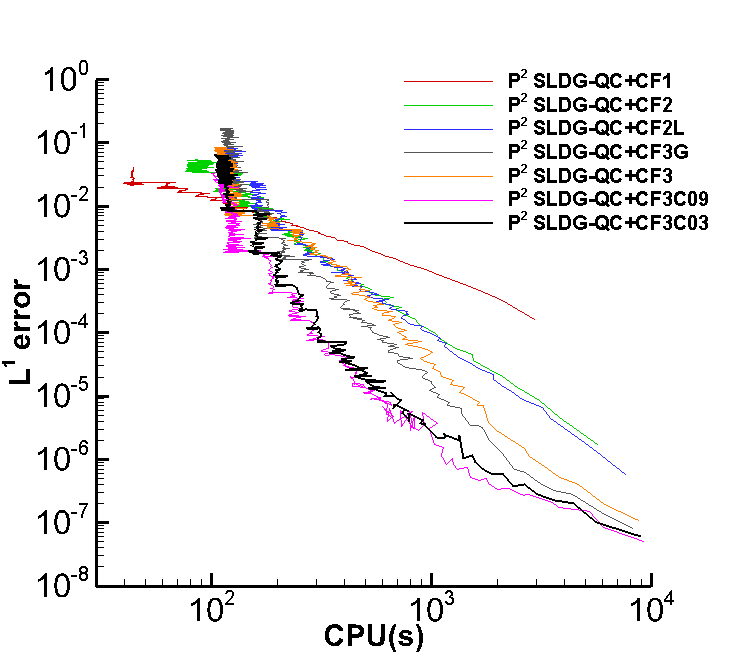}
\caption{  Plots of error versus the $CFL$ number (left) and  error versus CPU time (s) (right) for solving the Kelvin-Helmholtz instability at $T=5$.
Temporal order of convergence in $L^1$ norm of $P^2$ SLDG-QC+$P^3$ LDG  with  various temporal schemes (denoted by $P^2$ SLDG-QC+temporal scheme)  by comparing numerical solutions with
a reference solution from the corresponding scheme with $CFL = 0.1$. 
  The mesh of $120\times120$ is used.}
\label{KH_time}
\end{figure}

For this problem, the energy $\| \mathbf{E} \|_{L^2}^2 = \int_{\Omega} \mathbf{E}\cdot\mathbf{E}dxdy$ and enstrophy $\| \rho  \|_{L^2}^2 = \int_{\Omega} \rho^2 dxdy$ should remain constant in time.
Tracking relative deviations of these quantities provides a good measurement of the quality of numerical schemes.
We evaluate the performance of proposed schemes with a large $CFL=5$ for the Kelvin-Helmholtz instability for a long-time simulation.
Figure \ref{KH_surface} shows surface plots of the numerical solutions for the Kelvin-Helmholtz instability at $T=40$, and Figure \ref{KH_norms} reports the time evolution of the relative deviation of mass, energy and enstrophy of the numerical solutions.
From Figure \ref{KH_surface}, we can observe that $P^2$ SLDG-QC outperforms $P^1$ SLDG, with the same mesh.
From Figure \ref{KH_norms}, we can observe that
(1) the mass is conserved up to machine precision for each time step of the presented schemes;
(2)
higher order schemes performs better than lower order ones in terms of enstrophy; and different RKEI schemes have comparable performances.


\begin{figure}[h]
\centering                              
\includegraphics[height=45mm]{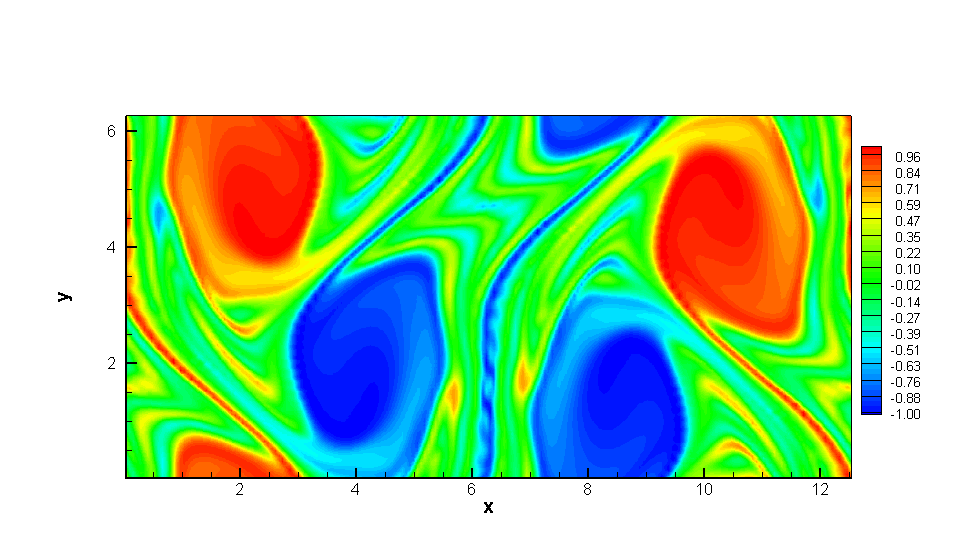}
\includegraphics[height=45mm]{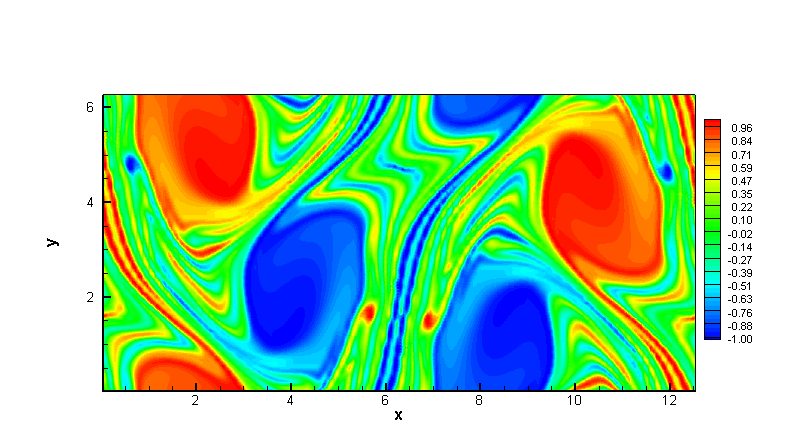}
\caption{Surface plots of the numerical solutions for the Kelvin-Helmholtz instability at $T=40$.
We use a mesh of $100\times100$ cells and $CFL=5$.
Left: $P^1$ SLDG+$P^2$ LDG+CF2.
Right: $P^2$ SLDG-QC+$P^3$ LDG+CF3C03.}
\label{KH_surface}
\end{figure}

%
%
%
%
%
%
%
%
\begin{figure}[h]
\centering                              
\includegraphics[height=60mm]{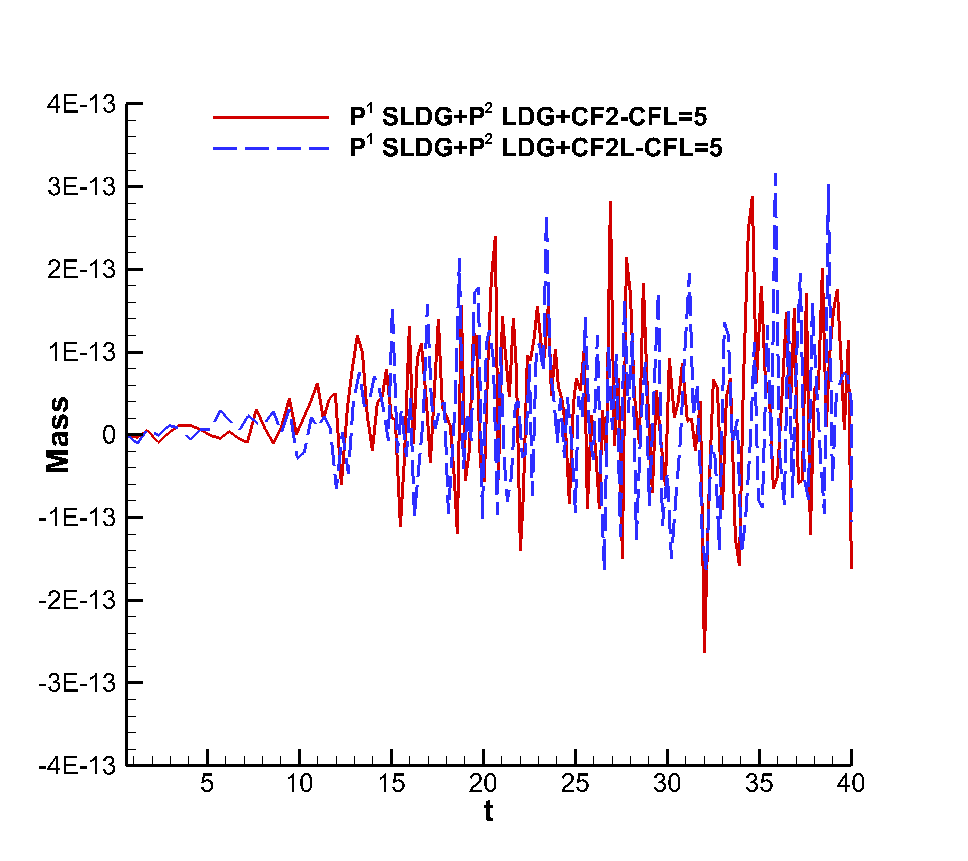}
\includegraphics[height=60mm]{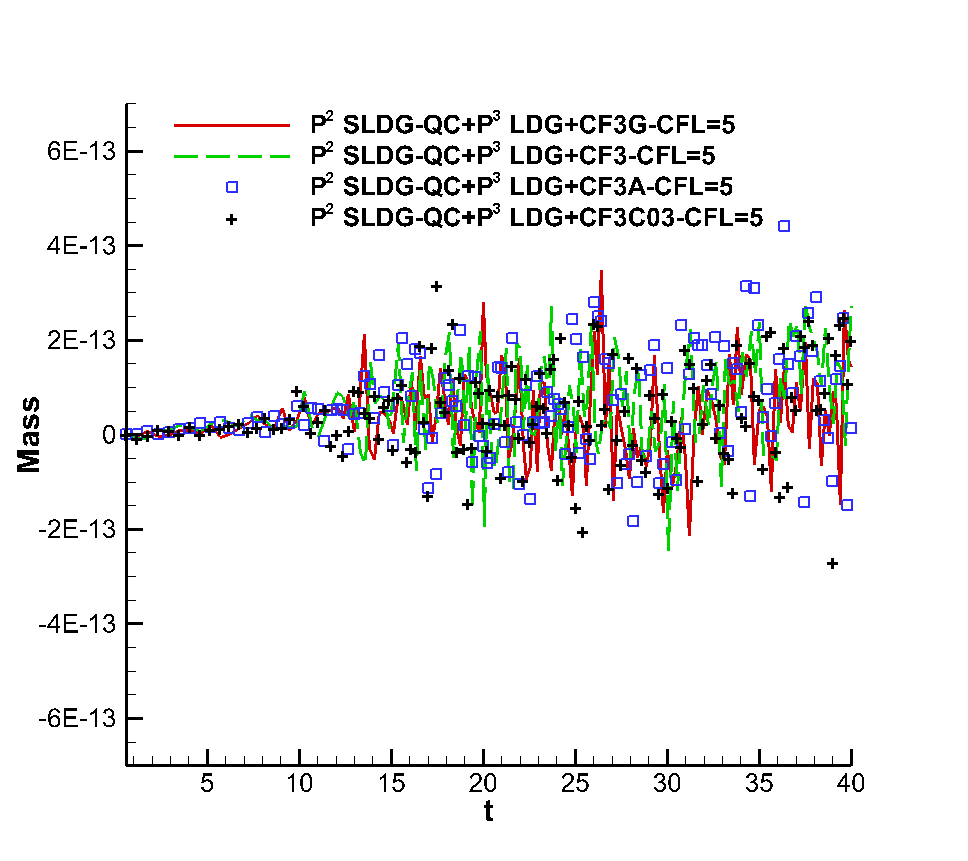}
\includegraphics[height=60mm]{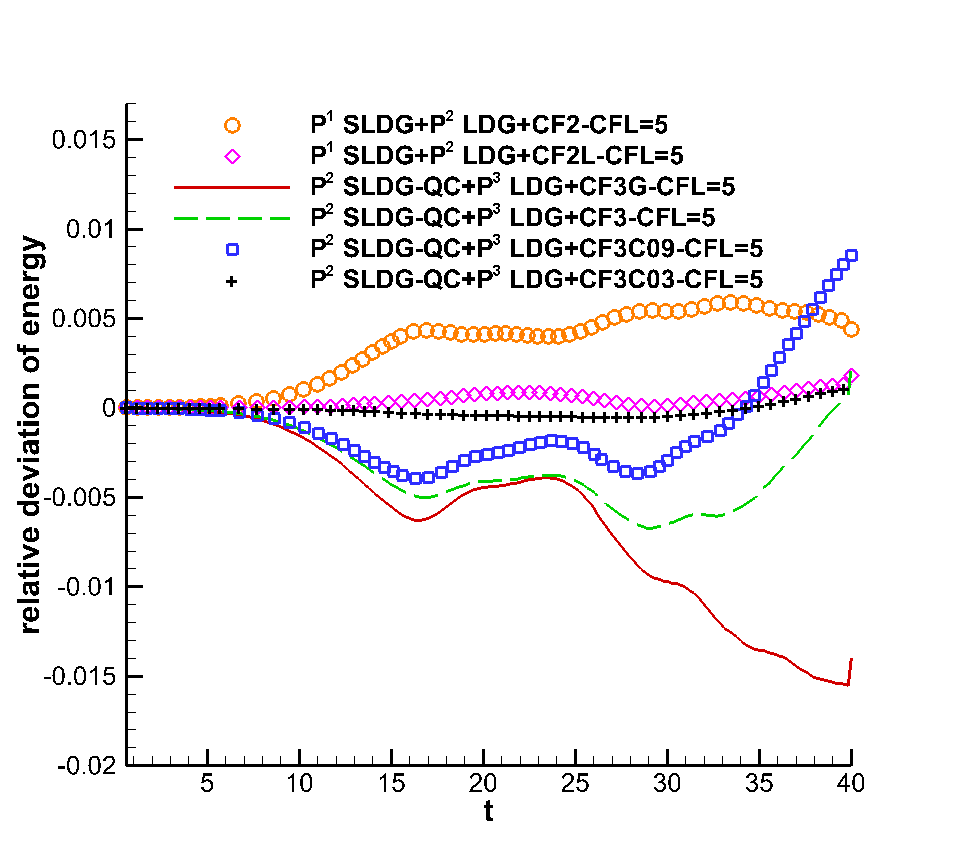}
\includegraphics[height=60mm]{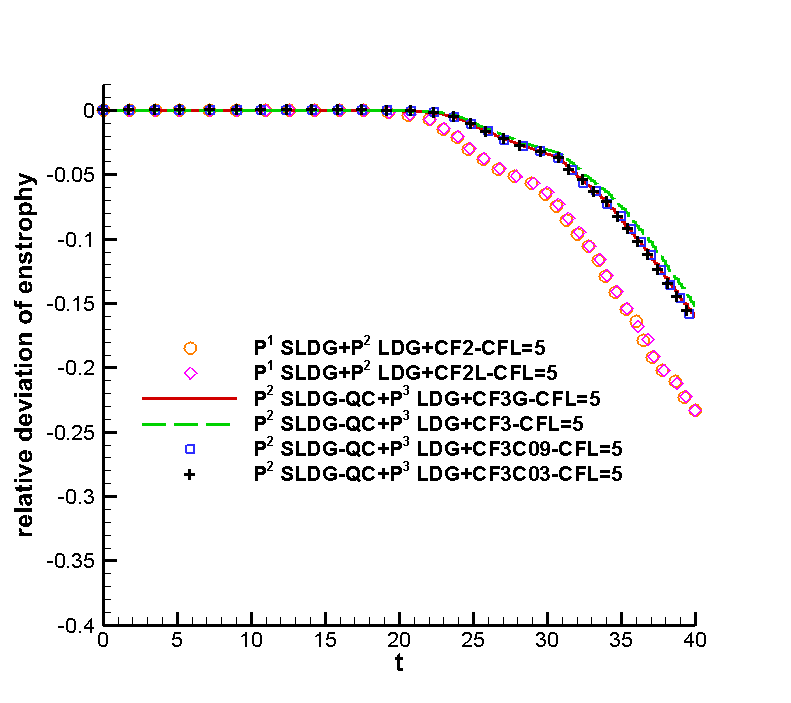}
\caption{Time evolutions of the relative deviation of mass (top left and right), energy (bottom left) and enstrophy (bottom right) for the proposed SLDG schemes for the Kelvin-Helmholtz instability problem. The mesh of $100\times100$ cells is used.}
\label{KH_norms}
\end{figure}
%
%

\begin{figure}[h]
\centering                              
\includegraphics[height=60mm]{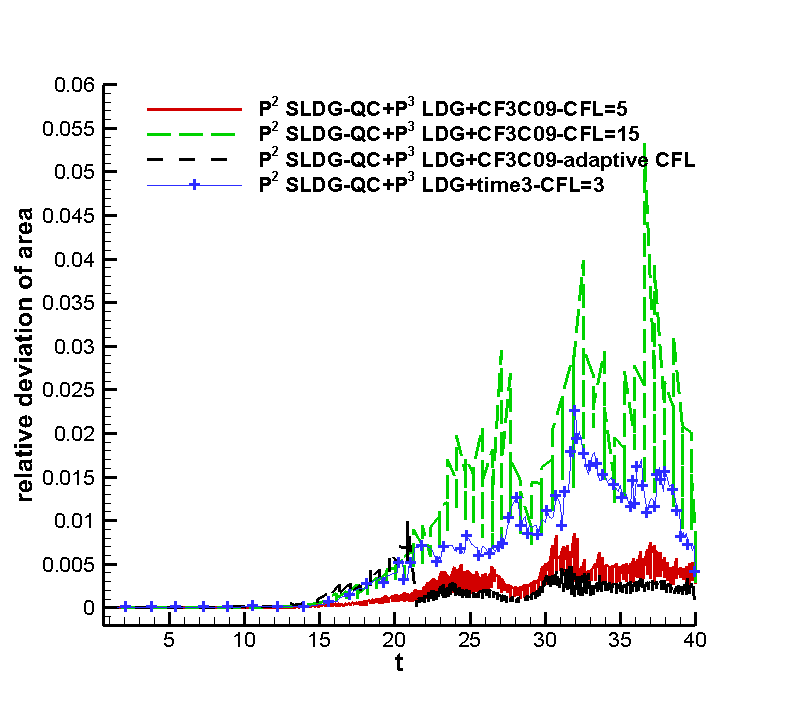}
\caption{Kelvin-Helmholtz instability. Use a mesh of $100\times100$ cells. Performances of $P^2$ SLDG-QC with different $CFL$s as well as the adaptive time-stepping algorithm in $L^\infty$ norm of the relative deviation of areas of upstream cells.}
\label{all_area}
\end{figure}

Next, we test the performance of proposed schemes  with the larger time-stepping size for the Kelvin-Helmholtz instability.
Note that the SLDG schemes with the very large time-stepping size might be subject to extreme distortion of upstream cells.
Due to the divergence-free constraint on the electric field of the guiding center Vlasov model, the areas of upstream cells should be preserved for the exact solution.
 If at the discrete level, the areas of upstream cells are preserved, the local maximum principle in terms of cell averages will be maintained; if   the  area  of  a  numerical  upstream  cell  greatly  deviates  from  the actual area, unphysical numerical oscillations may appear.
In particular, in Figure \ref{all_area}, by comparing  the relative deviation of area of $P^2$ SLDG-QC+$P^3$ LDG+time3 with $CFL=3$  and $P^2$ SLDG-QC+$P^3$ LDG+CF3C09 with the larger $CFL=5$, we observe that the latter one outperforms the former one.
It shows that the SLDG schemes with the third order exponential integrators can allow for a larger $CFL$.
Then we test the performance of proposed schemes  with a huge $CFL=15$.
We present surface plots of the numerical solutions of the schemes under $CFL=15$ at $T=40$ in Figure \ref{KH_huge}.
From Figure \ref{KH_huge}, we can find that   $P^2$ SLDG-QC+$P^3$ LDG+CF3C09 performs best, while $P^2$ SLDG-QC+$P^3$ LDG+CF3 performs worst in correctly resolving solution structure.
These observations show a good agreement with the accuracy comparison of these schemes presented in Figure~\ref{KH_time}.

\begin{figure}[h]
\centering                              
\includegraphics[height=45mm]{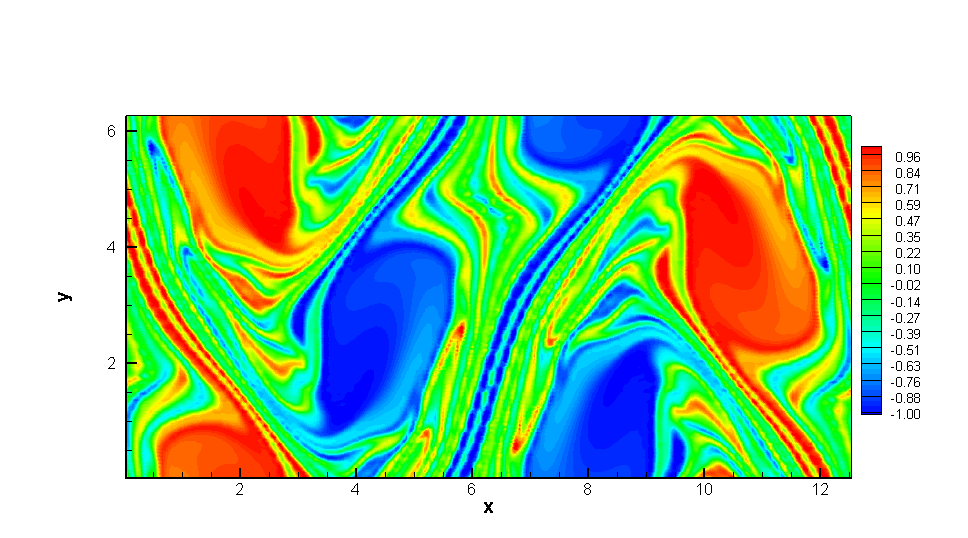}
\includegraphics[height=45mm]{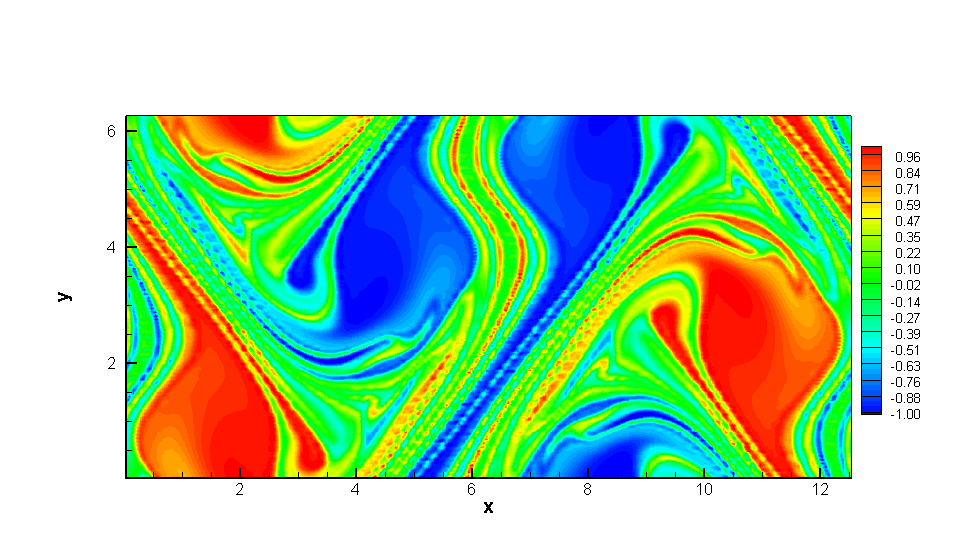}
\includegraphics[height=45mm]{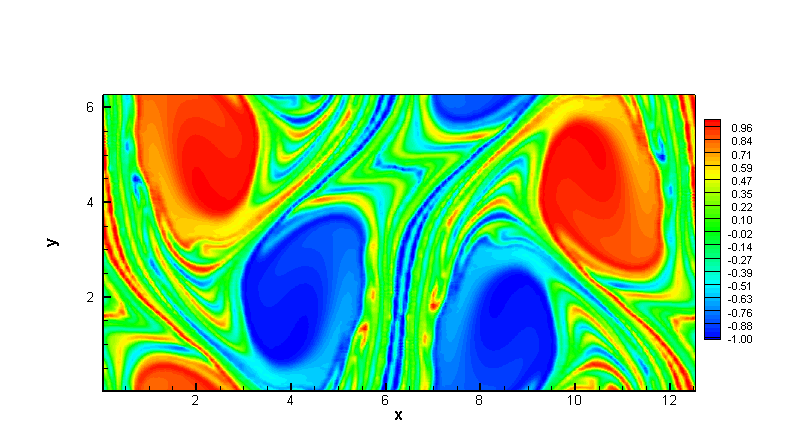}
\includegraphics[height=45mm]{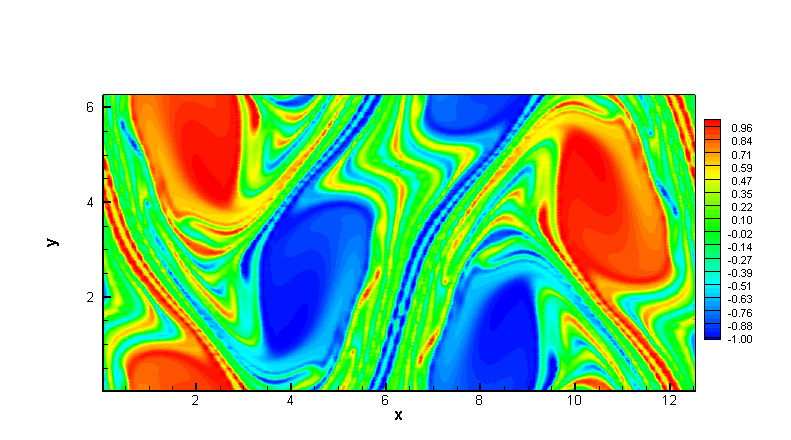}
\caption{ Surface plots of  numerical solutions for the Kelvin-Helmholtz instability at $T=40$.
We use a mesh of $100\times100$ cells and huge $CFL=15$.
Top left: $P^2$ SLDG-QC+$P^3$ LDG+CF3G. Top right: $P^2$ SLDG-QC+$P^3$ LDG+CF3.
Bottom left: $P^2$ SLDG-QC+$P^3$ LDG+CF3C09. Bottom right: $P^2$ SLDG-QC+$P^3$ LDG+CF3C03.}
\label{KH_huge}
\end{figure}


\noindent{\bf Numerical verification of the adaptive time-stepping algorithm for the nonlinear transport problems with the divergence-free constraint. }
It is observed that different RKEI schemes with $CFL=5$ have comparable performances, therefore we only present results from one second order and one third order RKEI schemes in Figure~\ref{KH_surface}; yet schemes with larger $CFL=15$ are observed to perform differently as shown in Figure \ref{KH_huge}. In Figure \ref{all_area}, we plot the schemes' performance in conserving the $L^\infty$ of the the upstream area under various settings (different schemes, CFLs, etc.). It is observed that,  compared to the scheme with $CFL=5$,
the area of upstream cells from the scheme with $CFL=15$ has larger deviation. We adopt the relative deviation of areas of upstream cells as an indicator for the adaptive time-stepping algorithm in Section \ref{section:adaptive}.

Now we test $P^2$ SLDG-QC+$P^3$ LDG+CF3C09 with the adaptive time-stepping algorithm.
In our numerical experience,  we find the choice of   $\delta_{M}=1\%$ to be optimal for this test case with  $T=40$.
We find the lower shreshold $\delta_{m}$ is less sensitive to the upper one.
In Figure  \ref{KH_adaptive}, we present the results of the scheme with $\delta_M = 1 $, and with two choices of
$\delta_{m}=0.1\%$ and  $\delta_{m}=0.05\%$.
 The left part of  Figure \ref{KH_adaptive} presents the 3D plots of the solution at $T=40$; and the right part of Figure \ref{KH_adaptive} plots of the adaptive $CFL$ versus time $T$, which showcases the effectiveness of the adaptive time stepping strategy when the solution gets more complex.
Figure \ref{norms_adaptive} shows the time evolutions of the relative deviation of energy and enstrophy the scheme with different settings.
The adaptive time-stepping scheme can conserve energy much better than the scheme with a fixed large $CFL$ number.

\begin{figure}[h]
\centering                              
\includegraphics[height=75mm]{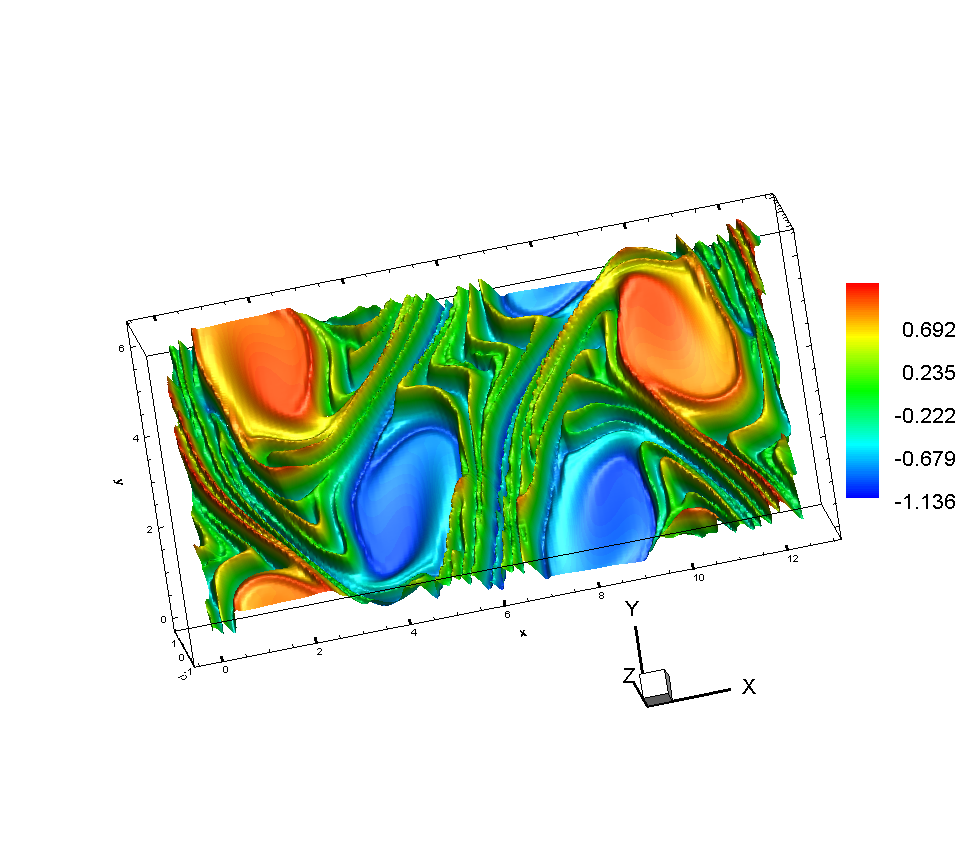}
\includegraphics[height=60mm]{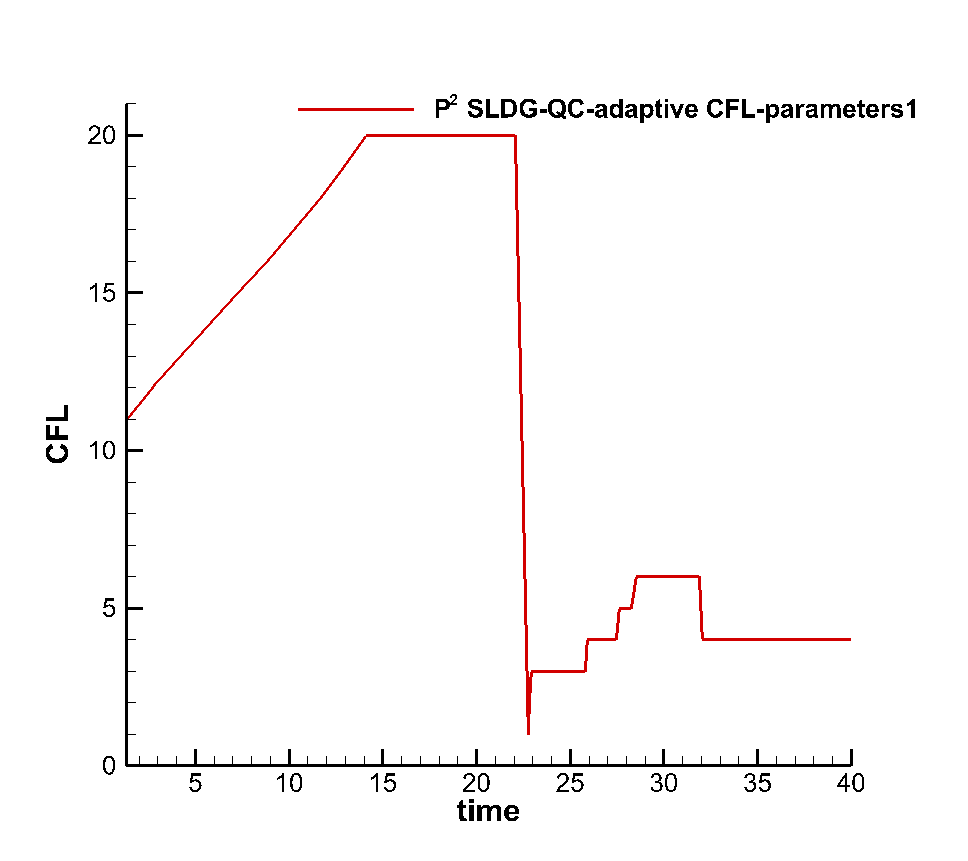}
\includegraphics[height=75mm]{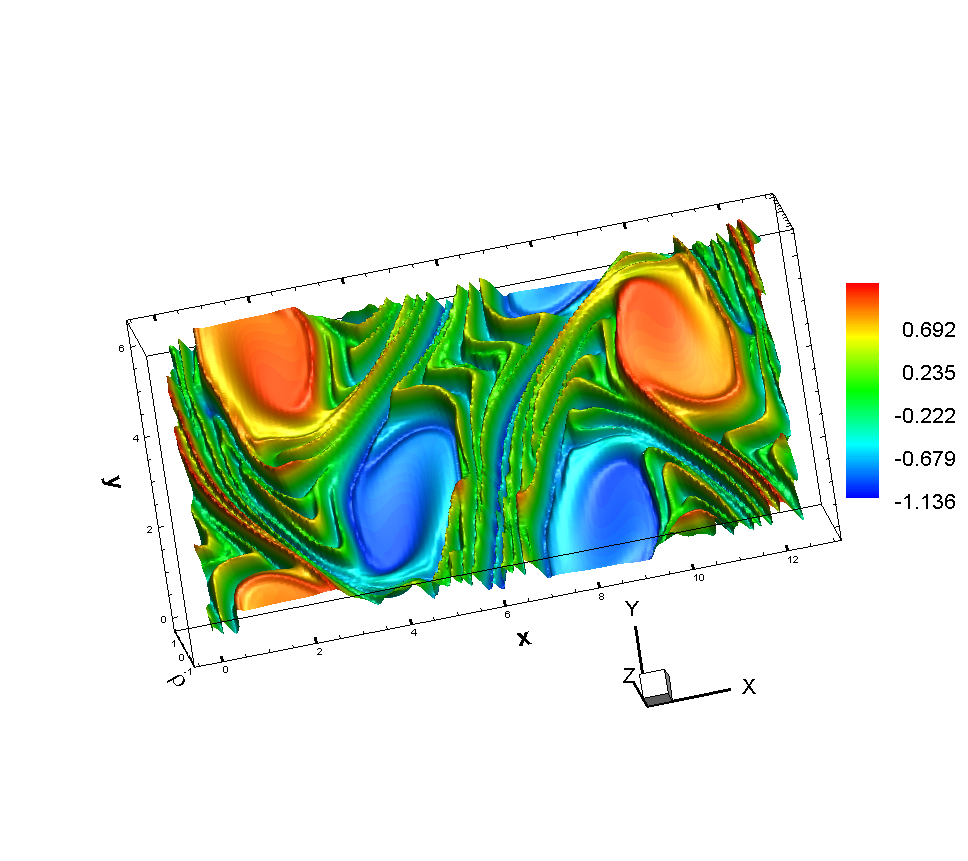}
\includegraphics[height=60mm]{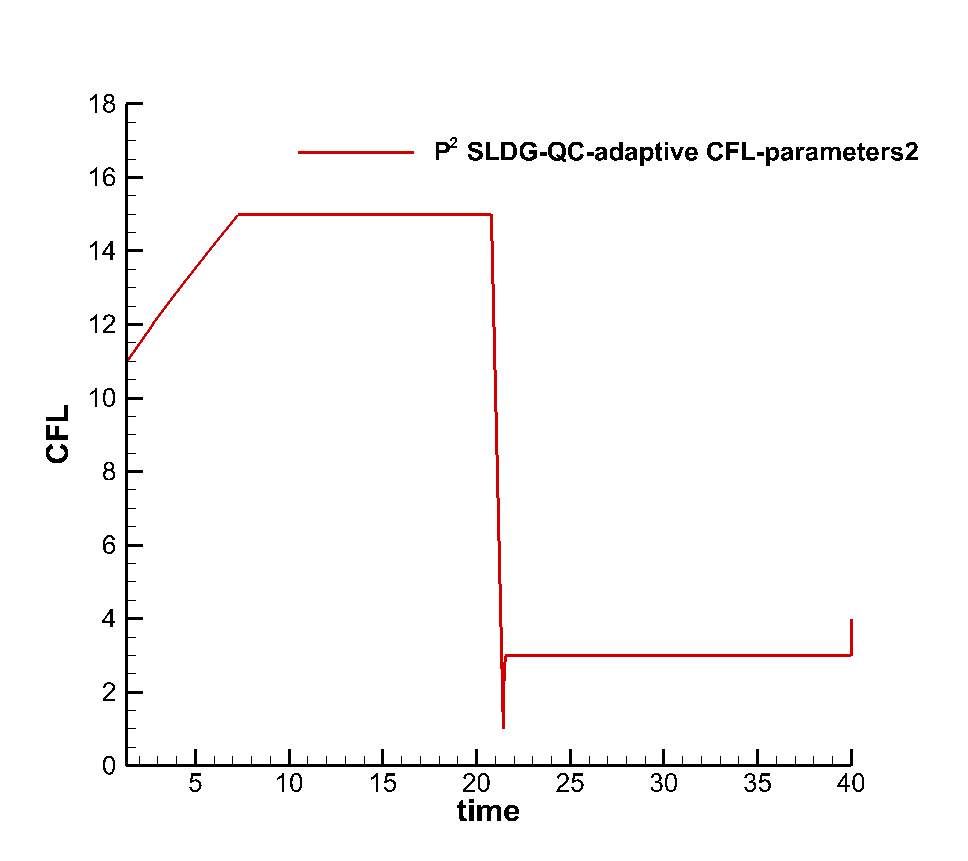}
\caption{
The numerical solution of $P^2$ SLDG-QC+$P^3$ LDG+CF3C09 with the adaptive time-stepping algorithm  for the Kelvin-Helmholtz instability at $T=40$, using a mesh of $100\times100$ cells.
Parameters 1: $\delta_{M}=1\%$, $\delta_{m}=0.1\%$;
Parameters 2: $\delta_{M}=1\%$, $\delta_{m}=0.05\%$.
Top/bottom left: the 3D plot of the solution of the SLDG using Parameters 1/Parameters 2.  Top/bottom right: the time evolution of $CFL$ versus time of the SLDG using Parameters 1/Parameters 2.
  }
\label{KH_adaptive}
\end{figure}

\begin{figure}[h]
\centering
\includegraphics[height=60mm]{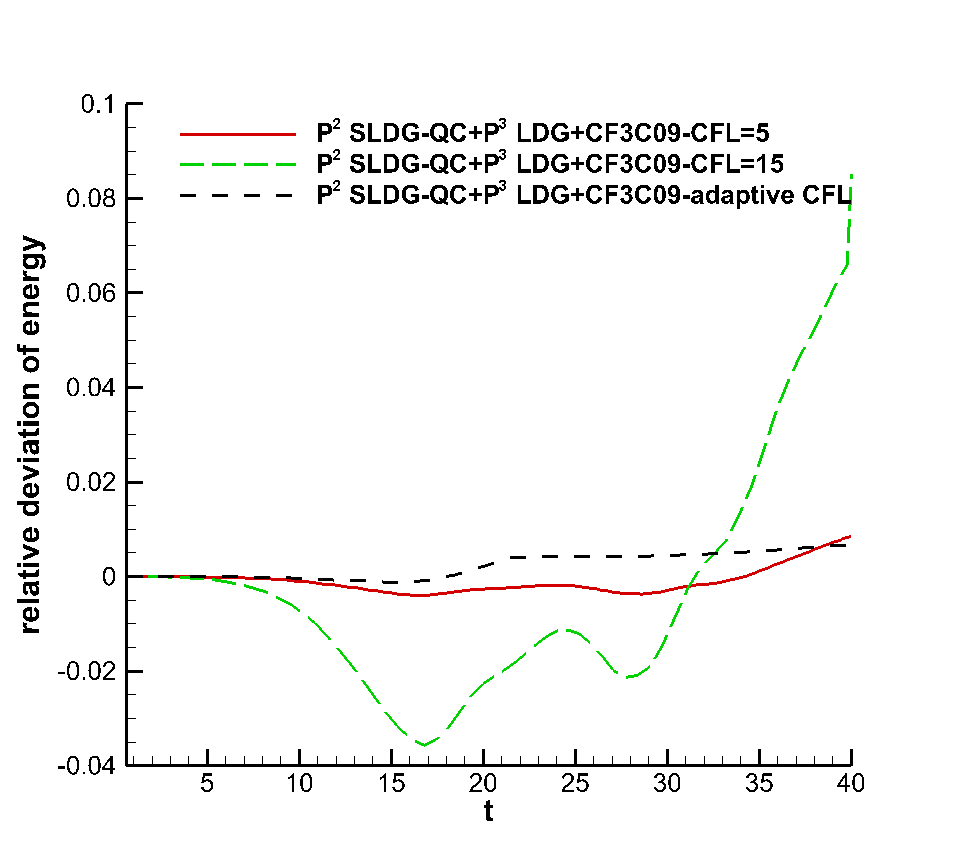}
\includegraphics[height=60mm]{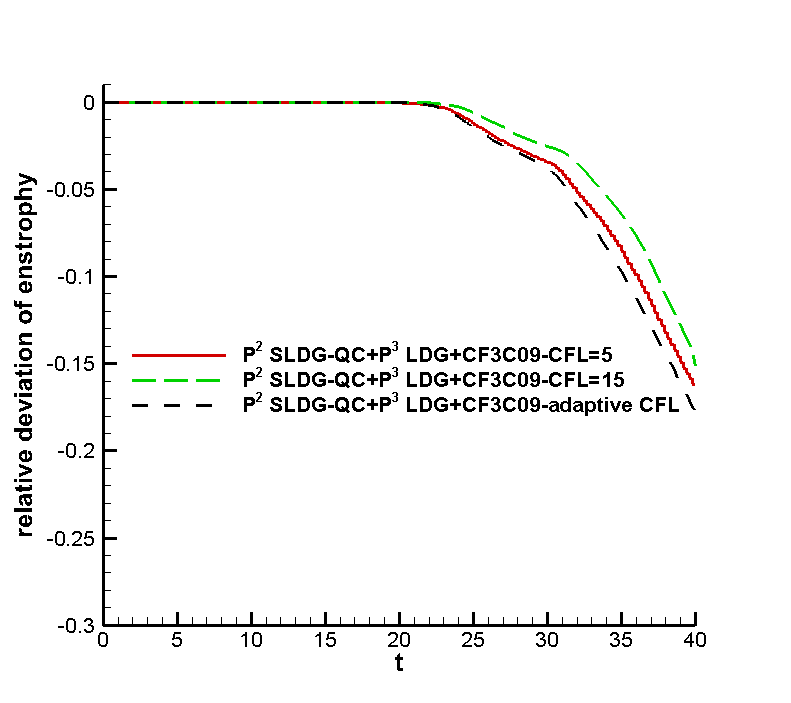}
\caption{The time evolutions of the relative deviation of energy (left) and enstrophy (right) for $P^2$ SLDG-QC+$P^3$ LDG+CF3C09 with the adaptive time-stepping algorithm for the Kelvin-Helmholtz instability problem, compared to that of the same scheme using a $CFL=5$ as well as that of the same scheme using a larger $CFL=15$.  }
\label{norms_adaptive}
\end{figure}

\end{exa}

\subsection{Preliminary numerical tests on Burgers' equation}
\label{sec4}

\begin{figure}[h!]
\centering
\includegraphics[height=70mm]{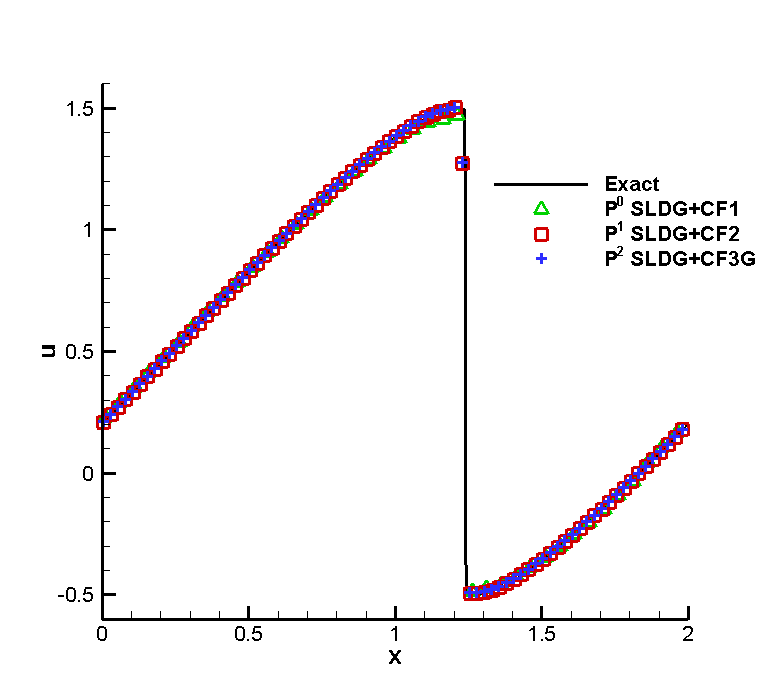}
\includegraphics[height=70mm]{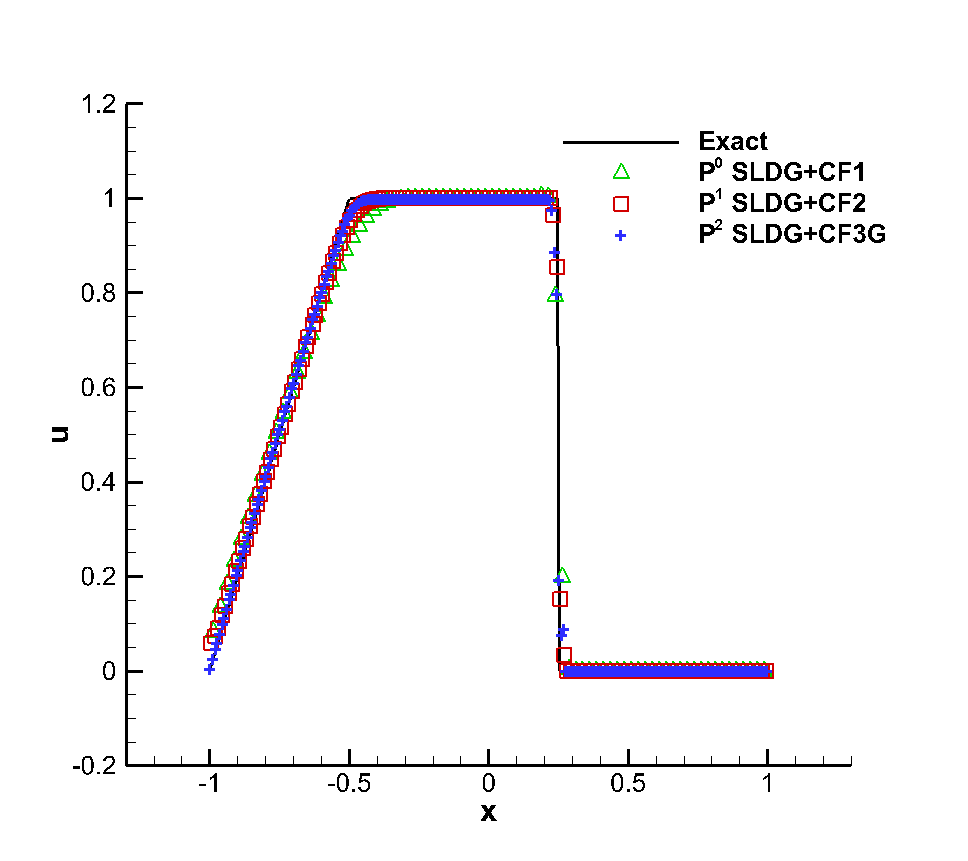}
\caption{Burgers's equation $u_t+(u^2/2)_x=0$.  Left: initial condition $u(x,0)=0.5+\sin(\pi x)$ at $t=1.5/\pi$. Right: a discontinuous initial condition \eqref{initial_13} and periodic boundary condition at
$t=0.5$.
The meshes of $80$ cells and $CFL=0.5$ are used.}
\label{burgers_shock}
\end{figure}

The focus of the current paper is on the nonlinear Vlasov dynamics, yet the SLDG-RKEI scheme can be applied to general nonlinear hyperbolic conservation laws.
In this subsection, we present our preliminary results on applying the SLDG-RKEI schemes to a simple Burgers' equation with shock developments.
Notice that for the nonlinear Vlasov dynamics, despite the development of filamentation structure as time evolves, shocks form from smooth initial data. We will show below that the development of shocks limit the time stepping size allowed for stability. We consider the 1D Burgers' equation,
\begin{equation}
u_t + \left(\frac12 u^2  \right)_x = 0,
\label{burgers1d}
\end{equation}
and rewrite it in the form of
\begin{equation}
u_t + \left(  P(u)   u \right)_x = 0.
\end{equation}
with $P(u) = 1/2 u$.
We first consider the initial condition $u(x,0)=0.5+\sin(\pi x)$. When $t=0.5/\pi$, the solution is still smooth. We present the errors and their corresponding orders of convergence in terms of $L^1$, $L^2$ and $L^\infty$ norms in Table \ref{burger1d_order}. Expected orders of convergence are observed. Notice that the presented $CFL$ for each scheme is the largest $CFL$ allowed for numerical stability, chopped off to one decimal place; and $\Delta t = \frac{ CFL \Delta x }{ \max_{u} f'(u) }$.  The maximum $CFL$ allowed for the Burgers' equation seems to be more limited than that for the nonlinear Vlasov dynamics; yet is still larger than those in a RKDG setting which is $1/(2k+1)$ with $k$ being the polynomial degrees.
The left panel of Figure \ref{burgers_shock} presents the solutions of the SLDG-RKEI schemes at $t=1.5/\pi$ after a shock develops using a mesh of $N=80$ with $CFL=0.5$.
We also consider a discontinuous initial condition,
\begin{equation}
u(x,0)
=
\begin{cases}
1   &  \text{if } -1\leq x<0,\\
0   &  \text{if } 0\leq x<1,
\end{cases}
\label{initial_13}
\end{equation}
with periodic boundary condition. The solution of this problem includes one shock and one rarefaction wave. The right panel of Figure \ref{burgers_shock} presents the numerical solutions that well capture the exact solution.
Note that the $CFL$ numbers taken here are ad hoc from numerical tests, and the stability property of the scheme for nonlinear problems is still subject to further investigation.

\begin{table}[!ht]
\caption{ Burgers's equation $u_t+(u^2/2)_x=0$ with the initial condition $u(x,0)=0.5+\sin(\pi x)$, $t=0.5/\pi$.
}
\vspace{0.1in}
\centering
\begin{tabular}{ l  cc cc cc }
\hline
Mesh  &{$L^1$ error} & Order    &{$L^2$ error} & Order   &{$L^\infty$ error} & Order  \\
 \hline
   \multicolumn{7}{l}{ $P^0$ SLDG+CF1 with $CFL=1.2$}
     \\

    40 &     2.75E-02 &      &     3.85E-02 &       &     1.72E-01 &     \\
    80 &     1.38E-02 &     0.99 &     1.96E-02 &     0.97 &     9.09E-02 &     0.92 \\
   160 &     6.96E-03 &     0.99 &     9.95E-03 &     0.98 &     4.71E-02 &     0.95 \\
   320 &     3.49E-03 &     1.00 &     5.02E-03 &     0.99 &     2.40E-02 &     0.97 \\
\hline

   \multicolumn{7}{l}{ $P^1$ SLDG+CF2 with $CFL=1.2$}
     \\

    40 &     1.81E-03 &       &     2.89E-03 &     &     1.66E-02 &       \\
    80 &     4.73E-04 &     1.93 &     7.57E-04 &     1.93 &     4.15E-03 &     2.00 \\
   160 &     1.20E-04 &     1.98 &     1.94E-04 &     1.97 &     1.02E-03 &     2.03 \\
   320 &     3.00E-05 &     2.00 &     4.90E-05 &     1.98 &     2.50E-04 &     2.03 \\
\hline
   \multicolumn{7}{l}{ $P^1$ SLDG+CF2L with $CFL=1.2$}
     \\

    40 &     1.79E-03 &       &     2.79E-03 &      &     1.15E-02 &     \\
    80 &     4.69E-04 &     1.94 &     7.37E-04 &     1.92 &     2.87E-03 &     2.00 \\
   160 &     1.19E-04 &     1.97 &     1.89E-04 &     1.96 &     7.00E-04 &     2.03 \\
   320 &     2.99E-05 &     1.99 &     4.78E-05 &     1.98 &     1.82E-04 &     1.95 \\
\hline
   \multicolumn{7}{l}{ $P^2$ SLDG+CF3G with $CFL=0.7$}
     \\

    40 &     1.11E-04 &      &     2.50E-04 &      &     2.14E-03 &       \\
    80 &     1.70E-05 &     2.71 &     4.55E-05 &     2.46 &     4.94E-04 &     2.12 \\
   160 &     2.59E-06 &     2.71 &     8.31E-06 &     2.45 &     1.27E-04 &     1.95 \\
   320 &     3.79E-07 &     2.77 &     1.50E-06 &     2.47 &     3.10E-05 &     2.04 \\
\hline
   \multicolumn{7}{l}{ $P^2$ SLDG+CF3  with $CFL=0.4$  }
   \\

    40 &     1.11E-04 &       &     2.42E-04 &      &     2.10E-03 &     \\
    80 &     1.68E-05 &     2.72 &     4.47E-05 &     2.44 &     4.89E-04 &     2.10 \\
   160 &     2.54E-06 &     2.73 &     8.20E-06 &     2.45 &     1.27E-04 &     1.95 \\
   320 &     3.76E-07 &     2.75 &     1.49E-06 &     2.46 &     3.09E-05 &     2.04 \\
\hline
      \multicolumn{7}{l}{ $P^2$ SLDG+CF3C09  with $CFL=0.7$}
     \\

    40 &     1.07E-04 &       &     2.40E-04 &     &     2.08E-03 &      \\
    80 &     1.65E-05 &     2.69 &     4.44E-05 &     2.43 &     4.88E-04 &     2.10 \\
   160 &     2.52E-06 &     2.72 &     8.18E-06 &     2.44 &     1.27E-04 &     1.95 \\
   320 &     3.73E-07 &     2.75 &     1.49E-06 &     2.46 &     3.09E-05 &     2.03 \\
\hline

   \multicolumn{7}{l}{ $P^2$ SLDG+CF3C03  with $CFL=0.7$}
     \\

    40 &     1.08E-04 &       &     2.41E-04 &       &     2.10E-03 &      \\
    80 &     1.65E-05 &     2.70 &     4.45E-05 &     2.44 &     4.89E-04 &     2.10 \\
   160 &     2.51E-06 &     2.72 &     8.20E-06 &     2.44 &     1.27E-04 &     1.95 \\
   320 &     3.72E-07 &     2.76 &     1.49E-06 &     2.46 &     3.09E-05 &     2.04 \\
\hline

\end{tabular}
\label{burger1d_order}
\end{table}

%


%
%
%
%
%

\section{Conclusion}\label{section:conclusion}

In this paper, we propose a high order SLDG-RKEI method for nonlinear Vlasov dynamics. Compared with previous work on semi Lagrangian methods, the new method could be systematically built up to be high order accurate in both spatial and temporal directions, mass conservative, computationally efficient in allowing extra large time stepping sizes and highly effective in resolving nonlinear Vlasov dynamics. Applications of the method to nonlinear Navier-Stokes system will be investigated in our future work.  
We also test the scheme on the nonlinear Burgers' equation and found CFL constraints similar to those from an Eulerian approach. Further study needs to be performed to better understand the stability of the method. 



\appendix
\section{Some commutator-free methods}
\label{append:a}

In this appendix we present some different Butcher tableaus of commutator-free methods from \cite{celledoni2003commutator,celledoni2009semi}, which are tested in this paper.
\begin{table}[htbp]
\centering
\begin{tabular}{r|c c } \label{CF2table}
$0$         &   &\\
$\frac12$	& $\frac12$  & 0\\\hline
            &  0         & 1
\end{tabular}

\caption{ CF2  }
\end{table}

\begin{table}[htbp]
\centering
\begin{tabular}{r|c c c} \label{CF2Ltable}
$0$         &           &            &  \\
$\gamma$	& $\gamma$  &            &  \\
     1      & $\delta$  & $1-\delta$ & \\\hline
            &   0       & $1-\gamma$ & $\gamma$
\end{tabular}

\caption{ CF2L. Here $\gamma=\frac{2-\sqrt{2} }{2}$ and $\delta=\frac{-2\sqrt{2} }{3}$.  }
\end{table}

%

\begin{table}[htbp]
\centering
\begin{tabular}{r|c c c} \label{CF3table}
$0$             &                &               &  \\
$\gamma$    	& $\gamma$      &               &  \\
   $1-\gamma$   &  $\gamma-1$    &  $2(1-\gamma)$            & \\ \hline
                &    0      & $\frac12-\phi$     & $\frac12+\phi$\\
                &    0      & $\phi$             & $  -\phi$
\end{tabular}

\caption{ CF3. Here $\gamma=\frac{3+\sqrt{3} }{6}$ and $\phi = \frac{1}{ 6(2\gamma-1) }$.  }
\end{table}

\begin{table}[htbp]
\centering
\begin{tabular}{r|c c c} \label{CF3Atable}
$0$             &                &               &  \\
$\gamma$    	& $\gamma$      &               &  \\
   $1-\gamma$   &  $\gamma-1$    &  $2(1-\gamma)$            & \\ \hline
                &   $\alpha$    & $\beta$     & $\sigma$\\
                &  $-\alpha$   & $\frac12-\beta$             & $\frac12-\sigma$
\end{tabular}

\caption{ CF3C09. Here $\sigma=(\alpha+\beta(1-2\gamma)-\frac13 )/(1-2\gamma),$ $\alpha=\frac12$, $\frac16$ and $\gamma=\frac{3+\sqrt{3} }{6}$.  }
\end{table}

\begin{table}[htbp]
\centering
\begin{tabular}{r|c c c} \label{CF3C03table}
$0$             &                &               &  \\
$\frac13$    	& $\frac13$      &               &  \\
   $\frac23$   &  $ 0$    &  $\frac23$            & \\ \hline
                &   $\frac13$    & $0$     & $0$\\
                &  $-\frac{1}{12}$   & $0$             & $\frac34$
\end{tabular}

\caption{ CF3C03. \cite{celledoni2003commutator}  }
\end{table}

\section{The quadratic-curved quadrilateral approximation}\label{append:b}

We highlight key components in the quadratic-curved quadrilateral
approximation to upstream cells

\begin{description}
  \item[Step 1.]

   \emph{Characteristics tracing.}  Locate the nine vertices of upstream element $A_{j}^\star$:
  $c_i^\star,\,i=1,\ldots,9$ by  tracking the characteristics backward to time $t^n$, i.e., solving the characteristics equations, for the nine vertices of $A_{j}$:
  $c_i,\,i=1,\ldots,9$  (see Figure \ref{schematic_2d_p2} (a)).

  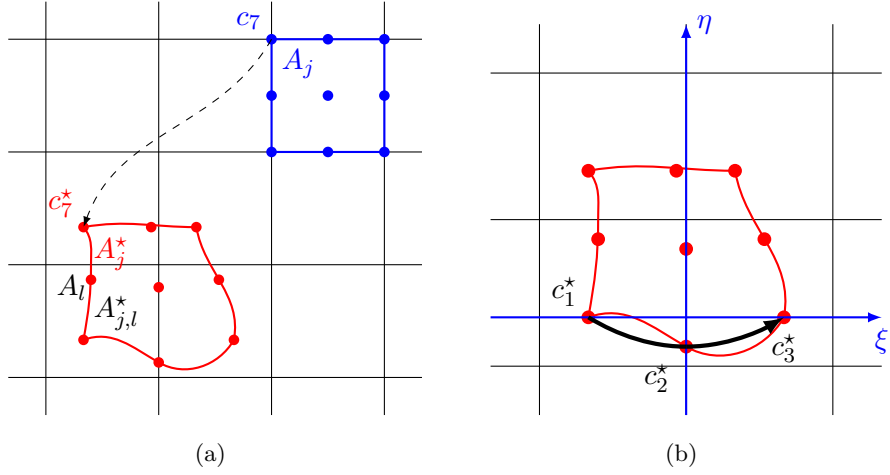
\begin{figure}[h]
\centering
\subfigure[]{
\begin{tikzpicture}
    \draw[black,thin] (0,0.5) node[left] {} -- (5.5,0.5)
                                        node[right]{};
    \draw[black,thin] (0,2.) node[left] {$$} -- (5.5,2)
                                        node[right]{};
    \draw[black,thin] (0,3.5) node[left] {$$} -- (5.5,3.5)
                                        node[right]{};
    \draw[black,thin] (0,5 ) node[left] {$$} -- (5.5,5)
                                        node[right]{};
    \draw[black,thin] (0.5,0) node[left] {} -- (0.5,5.5)
                                        node[right]{};
    \draw[black,thin] (2,0) node[left] {$$} -- (2,5.5)
                                        node[right]{};
    \draw[black,thin] (3.5,0) node[left] {$$} -- (3.5,5.5)
                                        node[right]{};
    \draw[black,thin] (5,0) node[left] {$$} -- (5,5.5)
                                        node[right]{};
    \fill [blue] (3.5,3.5) circle (2pt) node[] {};
    \fill [blue] (5,3.5) circle (2pt) node[] {};
    \fill [blue] (3.5,5) circle (2pt) node[below right] {$A_j$} node[above left] {$c_7$};
    \fill [blue] (5,5) circle (2pt) node[] {};

    \fill [blue] (3.5,4.25) circle (2pt) node[] {};
    \fill [blue] (5,4.25) circle (2pt) node[] {};
    \fill [blue] (4.25,4.25) circle (2pt) node[] {};
        \fill [blue] (4.25, 3.5) circle (2pt) node[] {};
    \fill [blue] (4.25,5) circle (2pt) node[] {};

     \draw[thick,blue] (3.5,3.5) node[left] {} -- (3.5,5)
                                        node[right]{};
      \draw[thick,blue] (3.5,3.5) node[left] {} -- (5,3.5)
                                        node[right]{};
       \draw[thick,blue] (3.5,5) node[left] {} -- (5,5)
                                        node[right]{};
        \draw[thick,blue] (5,3.5) node[left] {} -- (5,5)
                                        node[right]{};
    \fill [red] (1.,1) circle (2pt) node[above right,black] {$A_{j,l}^\star$};
    \fill [red] (3,1) circle (2pt) node[] {};
    \fill [red] (1,2.5) circle (2pt) node[below right] {$A_j^\star$} node[above left] {$c_7^\star$};
    \fill [red] (2.5,2.5) circle (2pt) node[] {};

     \draw[-latex,dashed](3.5,5)node[right,scale=1.0]{}
        to[out=240,in=70] (1,2.50) node[] {};

     \draw (0.5+0.01,2-0.01) node[fill=white,below right] {$A_l$};

     \draw [red,thick] (1,1)node[right,scale=1.0]{}
        to[out=20,in=150] (2,0.7) node[] {};

        \draw [red,thick] (2,0.7)node[right,scale=1.0]{}
        to[out=330,in=240] (3,1) node[] {};
             \draw [red,thick] (1,2.5)node[right,scale=1.0]{}
        to[out=310,in=90] (1.1,2) node[] {};
        \draw [red,thick] (1.1,2)node[right,scale=1.0]{}
        to[out=270,in=80] (1,1) node[] {};

        \draw [red,thick] (1,2.5)node[right,scale=1.0]{}
        to[out=10,in=180] (2.5,2.5) node[] {};

        \draw [red,thick] (3,1)node[right,scale=1.0]{}
        to[out=80,in=280] (2.5,2.5) node[] {};

       \fill [red] (2,0.7) circle (2pt) node[above right,black] {};
    \fill [red] (2,1.7) circle (2pt) node[] {};
    \fill [red] (1.9,2.5) circle (2pt) node[below right] {} node[above left] {};
    \fill [red] (1.1,1.8) circle (2pt) node[] {};
\fill [red] (2.8,1.8) circle (2pt) node[] {};

\end{tikzpicture}
}
\subfigure[]{

\begin{tikzpicture}[scale = 1.3]
    \draw[black,thin] (0,0.5) node[left] {} -- (4,0.5)
                                        node[right]{};
    \draw[black,thin] (0,2.) node[left] {$$} -- (4,2)
                                        node[right]{};
    \draw[black,thin] (0,3.5) node[left] {$$} -- (4,3.5)
                                        node[right]{};
    \draw[black,thin] (0.5,0) node[left] {} -- (0.5,4)
                                        node[right]{};
    \draw[black,thin] (2,0) node[left] {$$} -- (2,4)
                                        node[right]{};
    \draw[black,thin] (3.5,0) node[left] {$$} -- (3.5,4)
                                        node[right]{};

    \fill [red] (1.,1) circle (2pt) node[above right,black] {};
    \fill [red] (3,1) circle (2pt) node[] {};
    \fill [red] (1,2.5) circle (2pt) node[below right] {} node[above left] {};
    \fill [red] (2.5,2.5) circle (2pt) node[] {};

     \draw (0.5+0.01,2-0.01) node[fill=white,below right] {};

\draw [red,thick] (1,1)node[right,scale=1.0]{}
        to[out=20,in=150] (2,0.7) node[] {};

        \draw [red,thick] (2,0.7)node[right,scale=1.0]{}
        to[out=330,in=240] (3,1) node[] {};
             \draw [red,thick] (1,2.5)node[right,scale=1.0]{}
        to[out=310,in=90] (1.1,2) node[] {};
        \draw [red,thick] (1.1,2)node[right,scale=1.0]{}
        to[out=270,in=80] (1,1) node[] {};

        \draw [red,thick] (1,2.5)node[right,scale=1.0]{}
        to[out=10,in=180] (2.5,2.5) node[] {};

        \draw [red,thick] (3,1)node[right,scale=1.0]{}
        to[out=80,in=280] (2.5,2.5) node[] {};

               \fill [red] (2,0.7) circle (2pt) node[above right,black] {};
    \fill [red] (2,1.7) circle (2pt) node[] {};
    \fill [red] (1.9,2.5) circle (2pt) node[below right] {} node[above left] {};
    \fill [red] (1.1,1.8) circle (2pt) node[] {};
\fill [red] (2.8,1.8) circle (2pt) node[] {};

\draw[-latex,blue,thick](0,1)node[right,scale=1.0]{}
        to (4,1) node[below] {$\xi$};
\draw[-latex,blue,thick](2,0)node[right,scale=1.0]{}
        to (2,4) node[right] {$\eta$};

       \draw[black, ultra thick] (2,0.7) node[below left =2pt] {$c_2^\star$}
        parabola(1,1)node[above left,scale=1.0]{ $c_1^\star$ };
   \draw[-latex,black, ultra thick](2,0.7)node[above left,scale=1.0]{  }
         parabola (3,1) node[below =2pt] {$c_3^\star$};
\end{tikzpicture}

}
\caption{Schematic illustration of the SLDG formulation in two dimension. $P^2$ case.}
\label{schematic_2d_p2}
\end{figure}

  \item[Step 2.]

  \emph{Reconstructing faces of upstream elements.}
  Construct a quadratic curve to approximate each side of the upstream element. In particular, to construct the quadratic curve, $\wideparen{ c_1^\star,c_2^\star, c_3^\star }$ as shown in Figure \ref{schematic_2d_p2} (b), we do the following.

   \begin{description}
     \item[(a)]
     Construct  a coordinate transformation $x-y$ to $\xi-\eta$
      such that the coordinates of $c_1^\star$ and $c_3^\star$ are $(-1, 0)$ and $(1, 0)$ in $\xi-\eta$ space, respectively
      (see Figure \ref{schematic_2d_p2} (b)).
      Let $(x_1^\star, y_1^\star)$ and $(x_3^\star, y_3^\star)$ be the $x-y$ coordinate of $c_1^\star$ and $c_3^\star$, then
      the coordinate transformation is obtained as
     \begin{equation}
     \begin{cases}
     \xi(x,y) = a x + by + c,\\
     \eta(x,y) = bx -ay + d,
     \end{cases}
     \label{trans}
     \end{equation}
     where
     $$a=\frac{ 2( x_3^\star -x_1^\star  ) }{ (x_1^\star - x_3^\star  )^2 + (y_1^\star - y_3^\star  )^2  }, b=\frac{ 2( y_3^\star -y_1^\star  ) }{ (x_1^\star - x_3^\star  )^2 + (y_1^\star - y_3^\star  )^2  },$$
     $$c=\frac{   (x_1^\star)^2 - (x_3^\star)^2 + (y_1^\star)^2 - (y_3^\star)^2 }{ ( x_1^\star - x_3^\star  )^2 + (y_1^\star - y_3^\star  )^2  }, d=\frac{ 2(x_3^\star y_1^\star - x_1^\star y_3^\star  ) }{ (x_1^\star - x_3^\star  )^2 + (y_1^\star - y_3^\star  )^2  }.$$
     Its reverse transformation can be constructed accordingly:
     \begin{equation}
     \begin{cases}
     x = \frac{x_3^\star -x_1^\star }{2} \xi  + \frac{ y_3^\star -y_1^\star }{2} \eta + \frac{x_3^\star +x_1^\star }{2}, \\
     y = \frac{ y_3^\star -y_1^\star }{2} \xi  - \frac{ x_3^\star -x_1^\star }{2} \eta + \frac{ y_3^\star + y_1^\star }{2}.
     \end{cases}
     \label{reverse}
     \end{equation}

     \item[(b)]
    
     Get the $\xi-\eta$ coordinate for the point $c_2^\star$ as $(\xi_2, \eta_2)$. Based on $(-1,0)$, $(\xi_2,\eta_2)$ and $(1,0)$, we construct the parabola,
         \begin{equation}
         \wideparen{ c_1^\star,c_2^\star, c_3^\star }: \eta(x,y) = \frac{\eta_2}{ \xi_2^2-1 }(  \xi(x,y)^2-1  ).
         \label{parabola}
         \end{equation}
         Then
substitute \eqref{parabola}  into \eqref{reverse}, we have
     \begin{align}
     x(\xi) = \frac{x_3^\star -x_1^\star }{2} \xi  + \frac{ y_3^\star -y_1^\star }{2} \frac{\eta_2}{ \xi_2^2-1 }(  \xi^2-1  ) + \frac{x_3^\star +x_1^\star }{2}, \label{quadratic_x}\\
     y(\xi) = \frac{ y_3^\star -y_1^\star }{2} \xi  - \frac{ x_3^\star -x_1^\star }{2} \frac{\eta_2}{ \xi_2^2-1 }(  \xi^2-1  ) + \frac{ y_3^\star + y_1^\star }{2}. \label{quadratic_y}
     \end{align}
     
   \end{description}

  \item[Step 3.]
  
  \emph{Clipping algorithm.} As mentioned in \cite{cai2016high}, we perform a clipping algorithm for an upstream element by searching  its outer and inner segments, denoted as $ \mathcal{L}_q $ and $\mathcal{S}_q$.

  \item[Step 4.]

  \emph{Line integral evaluation.}
      The integral of inner line segments $ \sum_{q=1}^{N_i}
\int_{ \mathcal{S}_q } [Pdx +Q d y ]$ can be evaluated in the similar way as $P^1$ case.
The integral of out line segemtns $\sum_{q=1}^{N_o}
\int_{ \mathcal{L}_q } [P  dx +Q  d y ]  $ can be evaluated by the following parameterization on each side.
Assume that $\mathcal{L}_q$ is the part of the side $\wideparen{c_1^\star c_2^\star c_3^\star}$.
Hence,
\begin{equation}
\int_{  \mathcal{L}_q } [Pdx+ Qdy ]
= \int_{\xi^{(q)} }^{\xi^{(q+1)} }  [P( x(\xi,\eta),y(\xi,\eta) ) x'(\xi)+ Q( (\xi,\eta),y(\xi,\eta) ) y'(\xi) ]d\xi,
\end{equation}
where $(\xi^{(q)},\eta^{(q)})$ and $(\xi^{(q+1)},\eta^{(q+1)})$  are the starting point and the end point of $\mathcal{L}_q$ in $\xi-\eta$ coordinate, respectively.
The above integral can be done by the three point Gaussian quadrature.

\end{description}

\bibliographystyle{abbrv}
\bibliography{refer17}

\begin{thebibliography}{10}

\bibitem{arnold2002unified}
D.~Arnold, F.~Brezzi, B.~Cockburn, and L.~Marini.
\newblock {Unified analysis of discontinuous Galerkin methods for elliptic
  problems}.
\newblock {\em SIAM Journal on Numerical Analysis}, 39(5):1749--1779, 2002.

\bibitem{besse2017adaptive}
N.~Besse, E.~Deriaz, and {\'E}.~Madaule.
\newblock {Adaptive multiresolution semi-Lagrangian discontinuous Galerkin
  methods for the Vlasov equations}.
\newblock {\em Journal of Computational Physics}, 332:376--417, 2017.

\bibitem{besse2008wavelet}
N.~Besse, G.~Latu, A.~Ghizzo, E.~Sonnendr{\"u}cker, and P.~Bertrand.
\newblock {A wavelet-MRA-based adaptive semi-Lagrangian method for the
  relativistic Vlasov--Maxwell system}.
\newblock {\em Journal of Computational Physics}, 227(16):7889--7916, 2008.

\bibitem{bonaventura2018fully}
L.~Bonaventura, R.~Ferretti, and L.~Rocchi.
\newblock {A fully semi-Lagrangian discretization for the 2D incompressible
  Navier--Stokes equations in the vorticity-streamfunction formulation}.
\newblock {\em Applied Mathematics and Computation}, 323:132--144, 2018.

\bibitem{bosler2019conservative}
P.~A. Bosler, A.~M. Bradley, and M.~A. Taylor.
\newblock {Conservative Multimoment Transport along Characteristics for
  Discontinuous Galerkin Methods}.
\newblock {\em SIAM Journal on Scientific Computing}, 41(4):B870--B902, 2019.

\bibitem{butcher2008numerical}
J.~C. Butcher.
\newblock Numerical methods for ordinary differential equations, 2008.

\bibitem{cai2016high}
X.~Cai, W.~Guo, and J.-M. Qiu.
\newblock {A high order conservative semi-Lagrangian discontinuous Galerkin
  method for two-dimensional transport simulations}.
\newblock {\em Journal of Scientific Computing}, 73(2-3):514--542, 2017.

\bibitem{cai2018high}
X.~Cai, W.~Guo, and J.-M. Qiu.
\newblock {A high order semi-Lagrangian discontinuous Galerkin method for
  Vlasov-Poisson simulations without operator splitting}.
\newblock {\em Journal of Computational Physics}, 354:529--551, 2018.

\bibitem{cai2019high}
X.~Cai, W.~Guo, and J.-M. Qiu.
\newblock {A high order semi-Lagrangian discontinuous Galerkin method for the
  two-dimensional incompressible Euler equations and the guiding center Vlasov
  model without operator splitting}.
\newblock {\em Journal of Scientific Computing}, 79(2):1111--1134, 2019.

\bibitem{cai2019comparison}
X.~Cai, W.~Guo, and J.-M. Qiu.
\newblock {Comparison of semi-Lagrangian discontinuous Galerkin schemes for
  linear and nonlinear transport simulations}.
\newblock {\em Communications on Applied Mathematics and Computation},
  https://doi.org/10.1007/s42967-020-00088-0.

\bibitem{castillo2000priori}
P.~Castillo, B.~Cockburn, I.~Perugia, and D.~Sch{\"o}tzau.
\newblock {An a priori error analysis of the local discontinuous Galerkin
  method for elliptic problems}.
\newblock {\em SIAM Journal on Numerical Analysis}, 38(5):1676--1706, 2000.

\bibitem{celia1990eulerian}
M.~Celia, T.~Russell, I.~Herrera, and R.~Ewing.
\newblock {An Eulerian-Lagrangian localized adjoint method for the
  advection-diffusion equation}.
\newblock {\em Advances in Water Resources}, 13(4):187--206, 1990.

\bibitem{celledoni2009semi}
E.~Celledoni and B.~K. Kometa.
\newblock {Semi-Lagrangian Runge-Kutta exponential integrators for convection
  dominated problems}.
\newblock {\em Journal of Scientific Computing}, 41(1):139--164, 2009.

\bibitem{celledoni2016high}
E.~Celledoni, B.~K. Kometa, and O.~Verdier.
\newblock {High order semi-Lagrangian methods for the incompressible
  Navier--Stokes equations}.
\newblock {\em Journal of Scientific Computing}, 66(1):91--115, 2016.

\bibitem{celledoni2003commutator}
E.~Celledoni, A.~Marthinsen, and B.~Owren.
\newblock {Commutator-free Lie group methods}.
\newblock {\em Future Generation Computer Systems}, 19(3):341--352, 2003.

\bibitem{cheng2010preliminary}
A.~Cheng, K.~Wang, and H.~Wang.
\newblock {A preliminary study on multiscale ELLAM schemes for transient
  advection-diffusion equations}.
\newblock {\em Numerical Methods for Partial Differential Equations},
  26(6):1405--1419, 2010.

\bibitem{cockburn2003enhanced}
B.~Cockburn, M.~Luskin, C.-W. Shu, and E.~Suli.
\newblock {Enhanced accuracy by post-processing for finite element methods for
  hyperbolic equations}.
\newblock {\em Math. Comput.}, 72(242):577--606, 2003.

\bibitem{cockburn1998local}
B.~Cockburn and C.-W. Shu.
\newblock {The local discontinuous Galerkin method for time-dependent
  convection-diffusion systems}.
\newblock {\em SIAM Journal on Numerical Analysis}, 35(6):2440--2463, 1998.

\bibitem{cockburn2001runge}
B.~Cockburn and C.-W. Shu.
\newblock {Runge--Kutta discontinuous Galerkin methods for convection-dominated
  problems}.
\newblock {\em Journal of Scientific Computing}, 16(3):173--261, 2001.

\bibitem{crouseilles2014new}
N.~Crouseilles, P.~Glanc, S.~A. Hirstoaga, E.~Madaule, M.~Mehrenberger, and
  J.~P{\'e}tri.
\newblock {A new fully two-dimensional conservative semi-Lagrangian method:
  applications on polar grids, from diocotron instability to ITG turbulence}.
\newblock {\em The European Physical Journal D}, 68(9):252, 2014.

\bibitem{crouseilles2009conservative}
N.~Crouseilles, M.~Mehrenberger, and E.~Sonnendr{\"u}cker.
\newblock {Conservative semi-Lagrangian schemes for Vlasov equations}.
\newblock {\em Journal of Computational Physics}, 229(6):1927--1953, 2010.

\bibitem{crouseilles2011discontinuous}
N.~Crouseilles, M.~Mehrenberger, and F.~Vecil.
\newblock {Discontinuous Galerkin semi-Lagrangian method for Vlasov-Poisson}.
\newblock In {\em ESAIM: Proceedings}, volume~32, pages 211--230. EDP Sciences,
  2011.

\bibitem{degond2004modeling}
P.~Degond, L.~Pareschi, and G.~Russo.
\newblock {\em Modeling and computational methods for kinetic equations}.
\newblock Springer Science \& Business Media, 2004.

\bibitem{dimarco2015multiscale}
G.~Dimarco, R.~Loub{\`e}re, and V.~Rispoli.
\newblock {A multiscale fast semi-Lagrangian method for rarefied gas dynamics}.
\newblock {\em Journal of Computational Physics}, 291:99--119, 2015.

\bibitem{einkemmer2019performance}
L.~Einkemmer.
\newblock {A performance comparison of semi-Lagrangian discontinuous Galerkin
  and spline based Vlasov solvers in four dimensions}.
\newblock {\em Journal of Computational Physics}, 376:937--951, 2019.

\bibitem{einkemmer2014convergence}
L.~Einkemmer and A.~Ostermann.
\newblock Convergence analysis of a discontinuous galerkin/strang splitting
  approximation for the vlasov--poisson equations.
\newblock {\em SIAM Journal on Numerical Analysis}, 52(2):757--778, 2014.

\bibitem{grandgirard20165d}
V.~Grandgirard, J.~Abiteboul, J.~Bigot, T.~Cartier-Michaud, N.~Crouseilles,
  G.~Dif-Pradalier, C.~Ehrlacher, D.~Esteve, X.~Garbet, P.~Ghendrih, et~al.
\newblock {A 5D gyrokinetic full-f global semi-Lagrangian code for flux-driven
  ion turbulence simulations}.
\newblock {\em Computer Physics Communications}, 207:35--68, 2016.

\bibitem{groppi2016boundary}
M.~Groppi, G.~Russo, and G.~Stracquadanio.
\newblock {Boundary conditions for semi-Lagrangian methods for the BGK model}.
\newblock {\em Communications in Applied and Industrial Mathematics},
  7(3):138--164, 2016.

\bibitem{Guo2013discontinuous}
W.~Guo, R.~Nair, and J.-M. Qiu.
\newblock A conservative semi-{L}agrangian discontinuous {G}alerkin scheme on
  the cubed-sphere.
\newblock {\em Monthly Weather Review}, 142(1):457--475, 2013.

\bibitem{huang2016semi}
C.-S. Huang, T.~Arbogast, and C.-H. Hung.
\newblock {A semi-Lagrangian finite difference WENO scheme for scalar nonlinear
  conservation laws}.
\newblock {\em Journal of Computational Physics}, 322:559--585, 2016.

\bibitem{huang2012eulerian}
C.-S. Huang, T.~Arbogast, and J.~Qiu.
\newblock {An Eulerian--Lagrangian WENO finite volume scheme for advection
  problems}.
\newblock {\em Journal of Computational Physics}, 231(11):4028--4052, 2012.

\bibitem{kometa2011semi}
B.~K. Kometa.
\newblock {On Semi-Lagrangian Exponential Integrators and Discontinuous
  Galerkin Methods}.
\newblock In {\em AIP Conference Proceedings}, volume 1389, pages 1319--1322.
  AIP, 2011.

\bibitem{kormann2019massively}
K.~Kormann, K.~Reuter, and M.~Rampp.
\newblock {A massively parallel semi-Lagrangian solver for the six-dimensional
  Vlasov--Poisson equation}.
\newblock {\em The International Journal of High Performance Computing
  Applications}, 33(5):924--947, 2019.

\bibitem{lauritzen2010conservative}
P.~Lauritzen, R.~Nair, and P.~Ullrich.
\newblock {A conservative semi-Lagrangian multi-tracer transport scheme (CSLAM)
  on the cubed-sphere grid}.
\newblock {\em Journal of Computational Physics}, 229(5):1401--1424, 2010.

\bibitem{lee2016high}
D.~Lee, R.~Lowrie, M.~Petersen, T.~Ringler, and M.~Hecht.
\newblock {A high order characteristic discontinuous Galerkin scheme for
  advection on unstructured meshes}.
\newblock {\em J. Comput. Phys.}, 324:289--302, 2016.

\bibitem{li2016smoothness}
X.~Li, J.~K. Ryan, R.~M. Kirby, and C.~Vuik.
\newblock Smoothness-increasing accuracy-conserving (siac) filters for
  derivative approximations of discontinuous galerkin (dg) solutions over
  nonuniform meshes and near boundaries.
\newblock {\em Journal of Computational and Applied Mathematics}, 294:275--296,
  2016.

\bibitem{lin1997explicit}
S.~Lin and R.~Rood.
\newblock {An explicit flux-form semi-Lagrangian shallow-water model on the
  sphere}.
\newblock {\em Quarterly Journal of the Royal Meteorological Society},
  123(544):2477--2498, 1997.

\bibitem{peixoto2019semi}
P.~S. Peixoto and M.~Schreiber.
\newblock {Semi-Lagrangian Exponential Integration with application to the
  rotating shallow water equations}.
\newblock {\em SIAM Journal on Scientific Computing}, 41(5):B903--B928, 2019.

\bibitem{pironneau1982transport}
O.~Pironneau.
\newblock {On the transport-diffusion algorithm and its applications to the
  Navier-Stokes equations}.
\newblock {\em Numerische Mathematik}, 38(3):309--332, 1982.

\bibitem{qiu2017high}
J.-M. Qiu and G.~Russo.
\newblock A high order multi-dimensional characteristic tracing strategy for
  the {Vlasov--Poisson} system.
\newblock {\em Journal of Scientific Computing}, 71(1):414--434, 2017.

\bibitem{qiu_shu_sl}
J.-M. Qiu and C.-W. Shu.
\newblock {Conservative high order semi-Lagrangian finite difference WENO
  methods for advection in incompressible flow}.
\newblock {\em Journal of Computational Physics}, 230(4):863--889, 2011.

\bibitem{qiu2011positivity}
J.-M. Qiu and C.-W. Shu.
\newblock {Positivity preserving semi-Lagrangian discontinuous Galerkin
  formulation: Theoretical analysis and application to the Vlasov--Poisson
  system}.
\newblock {\em Journal of Computational Physics}, 230(23):8386--8409, 2011.

\bibitem{restelli2006semi}
M.~Restelli, L.~Bonaventura, and R.~Sacco.
\newblock {A semi-Lagrangian discontinuous Galerkin method for scalar advection
  by incompressible flows}.
\newblock {\em Journal of Computational Physics}, 216(1):195--215, 2006.

\bibitem{rossmanith2011positivity}
J.~A. Rossmanith and D.~C. Seal.
\newblock A positivity-preserving high-order semi-{L}agrangian discontinuous
  {G}alerkin scheme for the {V}lasov--{P}oisson equations.
\newblock {\em Journal of Computational Physics}, 230(16):6203--6232, 2011.

\bibitem{ryan2005extension}
J.~Ryan, C.-W. Shu, and H.~Atkins.
\newblock {Extension of a postprocessing technique for the discontinuous
  Galerkin method for hyperbolic equations with application to an aeroacoustic
  problem}.
\newblock {\em SIAM Journal on Scientific Computing}, 26(3):821--843, 2005.

\bibitem{shoucri1981two}
M.~M. Shoucri.
\newblock A two-level implicit scheme for the numerical solution of the
  linearized vorticity equation.
\newblock {\em International Journal for Numerical Methods in Engineering},
  17(10):1525--1538, 1981.

\bibitem{sonnendrucker2004vlasov}
E.~Sonnendr{\"u}cker, F.~Filbet, A.~Friedman, E.~Oudet, and J.-L. Vay.
\newblock Vlasov simulations of beams with a moving grid.
\newblock {\em Computer Physics Communications}, 164(1-3):390--395, 2004.

\bibitem{tumolo2012semi}
G.~Tumolo, L.~Bonaventura, and M.~Restelli.
\newblock {A semi-implicit, semi-Lagrangian, \emph{p}-adaptive discontinuous
  Galerkin method for the shallow water equations}.
\newblock {\em Journal of Computational Physics}, 232:46--67, 2013.

\bibitem{wang1999ellam}
H.~Wang, H.~Dahle, R.~Ewing, M.~Espedal, R.~Sharpley, and S.~Man.
\newblock {An ELLAM scheme for advection-diffusion equations in two
  dimensions}.
\newblock {\em SIAM Journal on Scientific Computing}, 20(6):2160--2194, 1999.

\bibitem{wang2006eulerian}
H.~Wang, W.~Zhao, R.~E. Ewing, M.~Al-Lawatia, M.~S. Espedal, and A.~S.
  Telyakovskiy.
\newblock {An Eulerian-Lagrangian solution technique for single-phase
  compositional flow in three-dimensional porous media}.
\newblock {\em Computers \& Mathematics with Applications}, 52(5):607--624,
  2006.

\bibitem{xiong2014high}
T.~Xiong, J.-M. Qiu, Z.~Xu, and A.~Christlieb.
\newblock {High order maximum principle preserving semi-Lagrangian finite
  difference WENO schemes for the Vlasov equation}.
\newblock {\em Journal of Computational Physics}, 273:618--639, 2014.

\bibitem{xiu2001semi}
D.~Xiu and G.~Karniadakis.
\newblock {A semi-Lagrangian high-order method for Navier-Stokes equations}.
\newblock {\em Journal of Computational Physics}, 172(2):658--684, 2001.

\bibitem{yang2020optimal}
Y.~Yang, X.~Cai, and J.-M. Qiu.
\newblock {Optimal convergence and superconvergence of semi-Lagrangian
  discontinuous Galerkin methods for linear convection equations in one space
  dimension}.
\newblock {\em Mathematics of Computation}, 89(325):2113--2139, 2020.

\bibitem{zhu2016h}
H.~Zhu, J.~Qiu, and J.-M. Qiu.
\newblock {An h-adaptive RKDG method for the Vlasov--Poisson system}.
\newblock {\em Journal of Scientific Computing}, 69(3):1346--1365, 2016.

\bibitem{zhu2017h}
H.~Zhu, J.~Qiu, and J.-M. Qiu.
\newblock {An h-Adaptive RKDG Method for the Two-Dimensional Incompressible
  {E}uler Equations and the Guiding Center {V}lasov Model}.
\newblock {\em Journal of Scientific Computing}, 73(2-3):1316--1337, 2017.

\end{thebibliography}

\end{document}